\input ./style/arxiv-general.cfg
\documentclass[MSNbibl,number,citesort,dvips]{arxbj}
\makeatletter
   \@ifpackageloaded{graphicx}{}{\usepackage{graphicx}}
\makeatother
\usepackage{upgreek}
\usepackage{accents}
\aid{0}
\volume{21}
\issue{3}
\pubyear{2015}
\firstpage{1719}
\lastpage{1759}
\doi{10.3150/14-BEJ619} % kopijuoti is 'New paper accepted'
\makeatletter
\renewcommand{\backslash}{\setminus}
\newcommand{\rrvert}{\vert}
\newcommand{\rrVert}{\Vert}
\newcommand{\llvert}{\vert}
\newcommand{\llVert}{\Vert}
\newcommand{\fracd}[2]{{(#1)}/{#2}}
\newcommand{\fraca}[2]{{#1}/{#2}}
\newcommand{\frace}[2]{{#1}/{(#2)}}
\renewcommand{\pi}{\uppi}
\newcommand{\eqref}[1]{(\ref{#1})}
\renewcommand{\mathring}{\accentset{\circ}}

\newproclaim{defi}{Definition}[section]
\newremark{rem}{Remark}[section]
\newtheorem{prop}{Proposition}[section]
\newproclaim{hyp}{Assumption}
\newremark{exple}{Example}[section]
\newtheorem{theorem}[prop]{Theorem}
\newtheorem{cor}[prop]{Corollary}
\newtheorem{lemma}[prop]{Lemma}

\newcommand{\EM}{}
\newcommand{\dC}{\EM{\mathbb{C}}}
\newcommand{\cD}{\EM{\mathcal{D}}}
\newcommand{\dE}{\EM{\mathbb{E}}}
\newcommand{\dN}{\EM{\mathbb{N}}}
\newcommand{\dP}{\EM{\mathbb{P}}}
\newcommand{\dR}{\EM{\mathbb{R}}}
\newcommand{\dZ}{\EM{\mathbb{Z}}}
\newcommand{\cB}{\EM{\mathcal{B}}}
\newcommand{\cC}{\EM{\mathcal{C}}}
\newcommand{\al}{\alpha}
\newcommand{\ga}{\gamma}
\newcommand{\be}{\beta}
\newcommand{\la}{\lambda}
\newcommand{\veps}{\varepsilon}
\makeatother

\begin{document}
\begin{frontmatter}

\title{Modulus of continuity of some conditionally sub-Gaussian
fields, application to stable random fields}
\runtitle{Modulus of continuity of some conditionally sub-Gaussian fields}

\begin{aug}
%%%% inicialai - be tarpu
\author[1]{\inits{H.}\fnms{Hermine}~\snm{Bierm\'e}\thanksref{1}\ead[label=e1]{Hermine.Bierme@mi.parisdescartes.fr}} \and
\author[2,3,4]{\inits{C.}\fnms{C\'eline}~\snm{Lacaux}\corref{}\thanksref{2,3,4}\ead[label=e2]{Celine.Lacaux@univ-lorraine.fr}}
%\author{\inits{}\fnms{}~\snm{}}
%%\runauthor{} %% auto
%\dedicated{}
\address[1]{MAP 5, CNRS UMR 8145, Universit\'e Paris Descartes, 45 rue
des Saints-P\`eres, 75006 Paris, France. \printead{e1}}
\address[2]{Institut \'Elie Cartan de Lorraine, Universit\'e de
Lorraine, UMR 7502, Vand\oe uvre-l\`es-Nancy, F-54506, France.
\printead{e2}}
\address[3]{CNRS, Institut \'Elie Cartan de Lorraine, UMR 7502,
Vand\oe uvre-l\`es-Nancy, F-54506, France}
\address[4]{INRIA, BIGS, Villers-l\`es-Nancy, F-54600, France}
\end{aug}

% HISTORY:
\received{\smonth{1} \syear{2013}}
\revised{\smonth{12} \syear{2013}}

% ABSTRACT
%
\begin{abstract}
In this paper, we study modulus of continuity and rate of convergence
of series of conditionally sub-Gaussian random fields. This framework
includes both classical series representations of Gaussian fields and
LePage series representations of stable fields. We enlighten their
anisotropic properties by using an adapted quasi-metric instead of the
classical Euclidean norm. We specify our assumptions in the case of
shot noise series where arrival times of a Poisson process are
involved. This allows us to state unified results for harmonizable
(multi)operator scaling stable random fields through their LePage
series representation, as well as to study sample path properties of
their multistable analogous.
\end{abstract}

% KEYWORDS
% visi is mazosios raides ir pagal abecele
%
\begin{keyword}
\kwd{H\"older regularity}
\kwd{operator scaling property}
\kwd{stable and multistable random fields}
\kwd{sub-Gaussian}
\end{keyword}
\end{frontmatter}

%s1 #&#
\section{Introduction}

In recent years, lots of new random fields
have been defined to propose new models for rough real data.
To cite a few of them, let us mention the (multi)fractional
Brownian fields (see, e.g., \cite{BEJARO}),
the linear and harmonizable (multi)fractional stable processes
\cite{StoevTaqqu04,Dozzi11} and some anisotropic fields such as the
(multi)fractional Brownian and stable sheets
\cite{AyacheLeger,AyacheRoueffXiao09} and the (multi)operator scaling
Gaussian and stable fields \cite{OSSRF,BLS11}.
In the Gaussian setting, sample path regularity relies on mean square
regularity. To study finer properties such as modulus of continuity, a
powerful technique consists in considering a representation of the
field as a series of random fields, using for instance Karhunen Loeve
decomposition (see \cite{Adler}, Chapter~3), Fourier or wavelet series
(as in \cite{Kahane60,AyacheXiao05}). This also allows generalizations
to non-Gaussian framework using for instance LePage series
\cite{LePage1,LePage2} for stable distributions (see, e.g., \cite{Kono}).
Actually, following previous works of LePage \cite{LePage2} and Marcus
and Pisier \cite{Marcus}, K\^ono and Maejima proved in \cite{Konot}
that, for $\alpha\in(0,2)$, an isotropic complex-valued $\alpha
$-stable random variable may be represented as a convergent shot noise
series of the form
%
%e1 #&#
%
\begin{equation}
\label{SN-Intro} \sum_{n=1}^{+\infty}T_n^{-1/\al}X_n,
\end{equation}
with $(T_n)_{n\ge1}$ the sequence of arrival times of a Poisson process
of intensity $1$, and $(X_n)_{n\ge1}$ a sequence of independent
identically distributed (i.i.d.) isotropic complex-valued random
variables, which is assumed to be independent of $(T_n)_{n\ge1}$ and
such that $\mathbb{E}(|X_1|^{\al})<+\infty$.
When $X_n=V_ng_n$ with $(g_n)_{n\ge1}$ a sequence of i.i.d.
Gaussian random variables independent of $(V_n,T_n)_{n\ge1}$, the
series may be considered as a conditional Gaussian series. This is one
of the main argument used in \cite{Kono,Dozzi11,BL09,BLS11} to study
the sample path regularity of some stable random fields. Another
classical representation consists in choosing $X_n=V_n\varepsilon_n$
with $(\varepsilon_n)_{n\ge1}$ a sequence of i.i.d. Rademacher random
variables that is, such that $\mathbb{P}(\varepsilon_n=1)=\mathbb
{P}(\varepsilon_n=-1)=1/2$. Both $g_n$ and $\varepsilon_n$ are
sub-Gaussian random variables.
Sub-Gaussian random variables have first been introduced in
\cite{Kahane60} for the study of random Fourier series. Their main property
is that their tail distributions behave like the Gaussian ones and then
sample path properties of sub-Gaussian fields
may be set as for Gaussian ones (see Theorem~12.16 of \cite{LedTal}, e.g.). In particular, they also rely on their mean square regularity.

In this paper, we study the sample path regularity of the
complex-valued series of conditionally sub-Gaussian fields defined as
%
%e2 #&#
%
\begin{equation}
\label{SGF-Intro} S(x)=\sum_{n=1}^{+\infty}
W_n(x) g_n, \qquad\mbox{for } x\in K_d
\subset \dR^d,
\end{equation}
with $\EM{{ (g_n )}}_{n\ge1}$ a sequence of independent symmetric
sub-Gaussian complex random variables, which is assumed independent of
$\EM{{ (W_n )}}_{n\ge1}$.
In this setting, we give sufficient assumptions on the sequence
$(W_n)_{n\ge1}$
to get an upper bound of the modulus of continuity of $S$ as well as a
uniform rate of convergence.
Then, we focus on shot noises series
\[
S\EM{{ (\al,u )}}=\sum_{n=1}^{+\infty}
T_n^{-1/\al} V_n(\al,u) g_n,
\qquad x=(\alpha,u)\in K_{d+1}\subset(0,2)\times{\dR}^d,
\]
with $\EM{{ (T_n )}}_{n\ge1}$
the sequence of arrival times of a Poisson process. Assuming the
independence of $\EM{{ (T_n )}}_{n\ge1}$, $\EM{{
(V_n )}}_{n\ge1}$ and $\EM{{ (g_n )}}_{n\ge1}$, we
state some more convenient conditions based on
moments of $V_n$ to ensure that the main assumptions of this paper are
fulfilled. In particular when $V_n(\alpha,u):=X_n$ is a symmetric
random variable, one of our main result gives a uniform rate of
convergence of the shot noise series \eqref{SN-Intro} in $\al$ on any
compact $K_1=[a,b]\subset(0,2)$, which improves the results obtained
in \cite{COLALE07} on the convergence of such series. In the framework
of LePage random series, which are particular examples of shot noise
series, we also establish that to improve the upper bound of the
modulus of continuity of $S$, one has the opportunity to use an other
series representation of $S$.
On the one hand, our framework allows to include in a general setting
some sample path regularity results already obtained in
\cite{BL09,BLS11} for harmonizable (multi)operator scaling stable random
fields. On the other hand, considering $\alpha$ as a function of $u
\in{\dR}^d$, we also investigate sample path properties of
multistable random fields that have been introduced in \cite{FLV09}.
To illustrate our results, we focus on harmonizable random fields.

The paper falls into the following parts. In Section~\ref{Preliminaries},
we recall definition and properties of sub-Gaussian
random variables and state our first assumption needed to ensure that
the random field $S$ is well-defined by \eqref{SGF-Intro}. We also
introduce a notion of anisotropic local regularity, which is obtained
by replacing\vadjust{\goodbreak} the isotropic Euclidean norm of ${\dR}^d$ by a
quasi-metric that can reveal the anisotropy of the random fields.
Section~\ref{MCGS} is devoted to our main results concerning both
local modulus of continuity of the random field $S$ defined by the
series \eqref{SGF-Intro} and rate of convergence of this series.
Section~\ref{SeSN} deals with the particular setting of shot noise
series, the
case of LePage series
being treated in Section~\ref{SeStable}. Then Section~\ref
{Applications} is devoted to the
study of the sample path regularity of stable or even multistable
random fields.
Technical proofs are postponed to \hyperref[Prel-A]{Appendix} for reader convenience.

%s2 #&#
\section{Preliminaries}\label{Preliminaries}

%s2.1 #&#
\subsection{Sub-Gaussian random variables}

Real-valued sub-Gaussian random variables have been defined by \cite
{Kahane60}. The structure of the class of these random variables and
some conditions for continuity of real-valued sub-Gaussian random
fields have been studied in \cite{SG81}. In this paper, we focus on
conditionally complex-valued sub-Gaussian random fields, where a
complex sub-Gaussian random variable is defined as follows.

%de2.1 #&#
%
\begin{defi}\label{SG} A complex-valued random variable $Z$ is
sub-Gaussian if there exists $s\in[0,+\infty)$ such that
%
%e3 #&#
%
\begin{equation}
\label{SGE} \forall z\in\dC, \qquad\dE\EM{{ \bigl(\mathrm{e}^{\Re\EM{{
(\overline {z}Z )}}}
\bigr)}}\le\mathrm{e}^{\fracd{s^2\EM
{{\llvert  z \rrvert }}^2}{2}}.
\end{equation}
\end{defi}

%re2.1 #&#
%
\begin{rem}
\label{comKa} This definition coincides also with complex sub-Gaussian
random variables as defined in \cite{Fukuda} in the more general
setting of random variables with values in a Banach space. Moreover,
for a real-valued random variable $Z$, it also coincides with the
definition in \cite{Kahane60}.
Kahane
\cite{Kahane60} called the smallest $s$ such that \eqref{SGE} holds
the Gaussian shift of the sub-Gaussian variable~$Z$. In this paper, if
\eqref{SGE} is fulfilled, we say that $Z$ is sub-Gaussian with
parameter $s$.
\end{rem}

%re2.2 #&#
%
\begin{rem} A complex-valued random variable $Z$ is sub-Gaussian if and
only if $\Re(Z)$ and $\Im(Z)$ are real sub-Gaussian random variables.
Note that if $Z$ is sub-Gaussian with parameter $s$ then $\mathbb
{E}(\Re(Z))=\mathbb{E}(\Im(Z))=0$ and $\mathbb{E}(\Re(Z)^2)\le
s^2$ as well as $\mathbb{E}(\Im(Z)^2)\le s^2$.
\end{rem}

The main property of sub-Gaussian random variables is that their tail
distributions decrease exponentially as the Gaussian ones (see Lemma~\ref{DecSG}). Moreover, considering convergent series of independent
symmetric sub-Gaussian random variables, a uniform
rate of decrease is available and the limit remains a sub-Gaussian
random variable. This result, stated below, is one of the main tool we
use to study sample path properties of conditionally sub-Gaussian
random fields.

%pr2.1 #&#
%
\begin{prop}
\label{SGP}
Let $\EM{{ (g_n )}}_{ n\ge1}$ be a sequence of independent
symmetric
sub-Gaussian random variables with parameter $s=1$. Let us consider a
complex-valued sequence $a=(a_n)_{n\ge1}$ such that
\[
\EM{{\llVert a\rrVert }}_{{\ell^2}}^2=\sum
_{n=1}^{+\infty} \EM {{\llvert a_n \rrvert
}}^2 <+\infty.
\]\vspace*{-\baselineskip}
\begin{enumerate}[2.]
\item[1.] Then, for any $t\in(0,+\infty)$,
$
\dP\EM{{ (\sup_{N\in\dN\backslash\EM{{ \{0 \}}}}
\EM{{\llvert  \sum_{n=1}^N a_n g_n \rrvert }}> t\EM{{\llVert  a\rrVert }}_{{\ell^2}} )}}\le8 \mathrm{e}^{- \fraca{t^2}{8}}$.
\item[2.] Moreover, the series $
\sum a_n g_n$ converges almost surely, and its limit $ \sum_{n=1}^{+\infty} a_n g_n$ is a sub-Gaussian random variable with parameter
$\EM{{\llVert  a\rrVert }}_{{\ell^2}}$.
\end{enumerate}
\end{prop}

\begin{pf} See Appendix~\ref{Prel-A}.
\end{pf}

%re2.3 #&#
%
\begin{rem} In the previous proposition, assuming that the parameter
$s=1$ is not restrictive since $a_n$ can be replaced by $a_ns_n$ and
$g_n$ by $g_n/s_n$ when $g_n$ is sub-Gaussian with parameter $s_n>0$.
\end{rem}
%
%s2.2 #&#
\subsection{Conditionally sub-Gaussian series}

In the whole paper, for $d\ge1$, $K_d=\prod_{j=1}^d[a_j,b_j]\subset
\dR^d$
is a compact $d$-dimensional interval and
for each integer $N\in\dN$, we consider
%
%e4 #&#
%
\begin{equation}
\label{SGSN} S_N(x)= \sum_{n=1}^{N}
W_n(x) g_n,\qquad x\in K_d,
\end{equation}
where $
\sum_{n=1}^0=0$ by convention and where
the sequence $(W_n,g_n)_{n\ge1}$ satisfies the following assumption.

%as1 #&#
%
\begin{hyp}
\label{HypTVg}
Let ${\EM{{ (g_n )}}}_{n\ge1}$ and
${\EM{{ (W_n )}}}_{n \ge1}$ be independent sequences of
random variables.
\begin{enumerate}[2.]
\item[1.] ${\EM{{ (g_n )}}}_{n \ge1}$ is a sequence of
independent symmetric
complex-valued sub-Gaussian random variables with parameter $s=1$.
\item[2.]${\EM{{ (W_n )}}}_{n \ge1}$ is a sequence of
complex-valued
continuous random fields defined on $K_d$ and such that
\[
\forall x\in K_d, \qquad\mbox{almost surely } \sum
_{n=1}^{+\infty} \EM{{\bigl\llvert W_n(x)
\bigr\rrvert }}^2<+\infty.
\]
\end{enumerate}
\end{hyp}

Under Assumption~\ref{HypTVg}, conditionally on $\EM{{
(W_n )}}_{n\ge1}$,
each $S_N$ is a sub-Gaussian random field defined on $K_d$. Moreover,
for each $x$, Proposition~\ref{SGP} and Fubini theorem lead to the
almost sure convergence of $S_N(x)$
as $N\to+\infty$.
The limit field $S$ defined by
%
%e5 #&#
%
\begin{equation}
\label{ChampS} S(x)=\sum_{n=1}^{+\infty}
W_n(x)g_n,\qquad x\in K_d\subset
\dR^d,
\end{equation}
is then a conditionally sub-Gaussian random field. In the sequel, we
study almost sure uniform convergence and rate of uniform convergence
of $\EM{{ (S_N )}}_{N\in\dN}$ as well as the sample path
properties of
$S$.

Assume first that each $g_n$ is a Gaussian random variable and that
each $W_n$ is a deterministic random field, which implies that $S$ is a
Gaussian centered random field. Then, it is well known that its sample
path properties are given by the behavior of
%
%e6 #&#
%
\begin{equation}
\label{VCS} s\EM{{ ({x,y} )}}:=\EM{{ \Biggl(\sum
_{n=1}^{+\infty} \EM {{\bigl\llvert W_n(x)-W_n(y)
\bigr\rrvert }}^2 \Biggr)}}^{1/2},\qquad x, y\in
K_d,
\end{equation}
since $s^2$ is
proportional to
the variogram
$
(x,y)\mapsto v(x,y):=\dE\EM{{ [\EM{{\llvert  S(x)-S(y) \rrvert }}^2 ]}}$.
In the following, we see that under Assumption~\ref{HypTVg}, the
behavior of $S$ is still linked with the behavior of the parameter $s$.
In this more general framework,
a key tool is to remark that conditionally on $\EM{{ (W_n
)}}_{n\ge1}$,
$S$ is a sub-Gaussian random field and the random variable
$S(x)-S(y)$ is sub-Gaussian with parameter $s(x,y)$.

We are particularly interested in anisotropic random fields $S$ (and
then anisotropic parameters $s$). Therefore, next section deals with an
anisotropic generalization of the classical H\"older regularity, that
is, with a notion of regularity which takes into account the anisotropy
of the fields under study.

%s2.3 #&#
\subsection{Anisotropic local regularity} \label{QM}

Let us first recall the notion of quasi-metric (see, e.g., \cite
{Quasimetric}), which is more adapted to our framework.

%de2.2 #&#
%
\begin{defi} A continuous function $\rho\dvtx{\mathbb R^d}\times
{\mathbb R^d}\rightarrow
[0,+\infty)$ is called a quasi-metric on ${\mathbb R^d}$ if
\begin{enumerate}[3.]
\item[1.]$\rho$ is faithful, that is, $\rho(x,y)=0$ iff $x=y$;
\item[2.]$\rho$ is symmetric, that is, $\rho(x,y)=\rho(y,x)$;
\item[3.]$\rho$ satisfies a quasi-triangle inequality: there exists a
constant $\kappa\ge1$ such that
\[
\forall x,y,z\in\dR^d,\qquad \rho(x,z)\le\kappa \bigl(\rho(x,y)+\rho (y,z)
\bigr).
\]
\end{enumerate}
\end{defi}

Observe that a continuous function $\rho$ is a metric on ${\mathbb
R^d}$ if and
only if $\rho$ is a quasi-metric on ${\mathbb R^d}$ which satisfies assertion
3 with $\kappa=1$.
In particular, the Euclidean distance
is an isotropic quasi-metric and its following anisotropic generalization
\[
(x,y)\mapsto\rho(x,y): = \Biggl(\sum_{i=1}^d|x_i-y_i|^{p/a_i}
\Biggr)^{1/p},\qquad\mbox{where } p>0 \mbox{ and } a_1,
\ldots,a_d>0,
\]
is also a quasi-metric.
Such quasi-metrics are particular cases of the following general example.

%ex2.1 #&#
%
\begin{exple}
\label{tauE}
Let us consider $E$ a real $d\times d$ matrix whose eigenvalues have
positive real parts and
$\tau_{ {E}}\dvtx{\mathbb R^d}\rightarrow\dR^+$ a continuous even
function such that
\begin{enumerate}[(ii)]
\item[(i)] for all $x\neq0$, $\tau_{ {E}}(x)>0$;
\item[(ii)] for all $r>0$ and all $x\in{\mathbb R^d}$, $\tau_{
{E}}(r^Ex)=r\tau_{ {E}}(x)$ with $r^E =\exp\EM{{ ((\ln r)
E )}}$.\vadjust{\goodbreak}
\end{enumerate}
The classical example of such a function is the radial part of polar
coordinates with respect to $E$ introduced in Chapter~6 of \cite
{thebook}. Other examples have been given in
\cite{OSSRF}.

Let us consider the continuous function $\rho_{ {E}}$, defined on
${\mathbb R^d}\times{\mathbb R^d}$ by
\[
\rho_{ {E}}(x,y)=\tau_{ {E}}(x-y).
\]
Then, by definition of $\tau_{{E}}$, $\rho_{ {E}}$ is faithful and
symmetric. Moreover, by Lemma~2.2 of \cite{OSSRF}, $\rho_{ {E}}$
also satisfies a quasi-triangle inequality. Hence, $\rho_{ {E}}$
is a quasi-metric on ${\mathbb R^d}$ and it is adapted to study
operator scaling
random fields (see, \cite{OSSRF,BL09}, e.g.).
\end{exple}

Let us remark that
since $\rho_{ {E}}^\beta$ defines a quasi-metric for $E/\beta$
whatever $\beta>0$ is, we may restrict our study to matrix $E$ whose
eigenvalues have real parts greater than one. Then, by Proposition~3.5
of \cite{BLS11},
there exist $0<\underline{H}\le\overline{H}\le1$ and two
constants
$c_{{2,1}},c_{{2,2}}\in(0,\infty)$
such that for all $x,y \in{\mathbb R^d}$,
\[
c_{{2,1}}\min\bigl(\|x-y\|^{\overline{H}},\|x-y\|^{\underline{H}}\bigr)
\le \rho_{ {E}}(x,y)\le c_{{2,2}}\max\bigl(\|x-y
\|^{\overline{H}},\|x-y\| ^{\underline{H}}\bigr),
\]
where $\|\cdot\|$ is the Euclidean norm on ${\mathbb R^d}$.
In \cite{BLS11,BL09}, this comparison is one of the main tool in the
study of the regularity of some stable anisotropic random fields.
Therefore, throughout the paper, we consider a quasi-metric $\rho$
such that there exist $0<\underline{H}\le\overline{H}\le1$ and two
constants
$c_{{2,1}},c_{{2,2}}\in(0,\infty)$ such that for all $x,y \in
{\mathbb R^d}
$, with $\|x-y\|\le1$,
%
%e7 #&#
%
\begin{equation}
\label{controlrho} c_{{2,1}}\|x-y\|^{\overline{H}}\le\rho(x,y)\le
c_{{2,2}}\|x-y\| ^{\underline{H}}.
\end{equation}
Before we introduce the anisotropic regularity used in the following,
let us briefly comment this assumption.

%re2.4 #&#
%
\begin{rem}
\begin{enumerate}[3.]
\item
The upper bound is needed in the sequel to construct a particular
$2^{-k}$ net for $\rho$,
whose cardinality
can be estimated using the lower bound.
\item Using the quasi-triangle inequality satisfied by $\rho$ and its
continuity, one deduces from \eqref{controlrho} that for any non-empty
compact set $K_d=\prod_{j=1}^d [a_j,b_j]\subset\dR^d$, there exist
two finite positive constants $c_{{2,1}}(K_d)$ and $c_{{2,2}}(K_d)$
such that for all $x,y \in K_d$,
%
%e8 #&#
%
\begin{equation}
\label{controlrho2} c_{{2,1}}(K_d)\|x-y\|^{\overline{H}}\le
\rho(x,y)\le c_{{2,2}}(K_d)\|x-y\|^{\underline{H}}.
\end{equation}
\item It is not restrictive to assume that $\overline{H}\le1$ since
for any $c>0$, $\rho^c$ is also a quasi-metric.
\end{enumerate}
\end{rem}

We will consider the following anisotropic local and uniform regularity
property.

%de2.3 #&#
%
\begin{defi} Let $\be\in(0,1]$ and $\eta\in\dR$. Let $x_0\in K_d$
with $K_d\subset\dR^d$. A real-valued function $f$ defined on $K_d$
belongs to ${\mathcal H}_{\rho,K_d}(x_0,\beta, \eta)$
if there exist $\gamma\in(0,1)$ and
$C\in(0,+\infty)$ such that
\[
\bigl |f\EM{{ (x )}}-f\EM{{ (y )}}\bigr |\le C \rho (x,y)^{\beta}\bigl |\log\bigl(\rho(x,y)
\bigr)\bigr |^{\eta}
\]
for all $x,y\in B(x_0,\gamma)\cap K_d=\{z\in K_d; \|z-x_0\|\le\gamma
\}$.
Moreover $f$ belongs to
${\mathcal H}_{\rho}(K_d,\beta, \eta)$ if there exists $C\in
(0,+\infty)$ such that
\[
\forall x,y\in K_d,\qquad\bigl |f\EM{{ (x )}}-f\EM{{ (y )}}\bigr |\le C
\rho(x,y)^{\beta}\EM{{ \bigl[\log\bigl(1+\rho (x,y)^{-1}\bigr)
\bigr]}}^{\eta}.
\]
\end{defi}

%re2.5 #&#
%
\begin{rem}\label{AR1}
\begin{enumerate}[3.]
\item If $f \in{\mathcal H}_{\rho,K_d}(x_0,\beta, \eta)$, then $f$
is continuous at $x_0$.
Moreover, since $h^\beta\EM [\log(1+  h^{-1}) ]^{\eta
}\sim_{h\to
0_+} h^\beta\EM{{\llvert  \log(h) \rrvert }}^{\eta}$ and since $\rho
$ satisfies
equation \eqref{controlrho}, $f \in{\mathcal H}_{\rho,K_d}(x_0,\beta
, \eta)$
if and only if for some $\gamma>0$, $f \in{\mathcal H}_{\rho
}(B(x_0,\gamma)\cap K_d,\beta, \eta)$.
\item If $f \in{\mathcal H}_{\rho}(K_d,\beta, \eta)$, then $f\in
{\mathcal H}_{\rho,K_d}(x_0,\beta, \eta)$ for all $x_0\in K_d$. The
converse is also true since $K_d$ is a compact. This follows from the
Lebesgue's number lemma and the boundedness of the continuous function
$f$ on the compact set $K_d$ (see Lemma~\ref{LemReco} stated in the
\hyperref[appB]{Appendix} for an idea of the proof).
\item A function in ${\mathcal H}_{\rho}(K_d,\beta, 0)$ may be view
as a Lipschitz function on an homogeneous space \cite{Macias79}. Note
also that
when $\rho$ is the Euclidean distance, for any $\beta\le1$ and $\eta
\le0$, the set ${\mathcal H}_{\rho}({K_d},\beta, \eta)$
(resp., ${\mathcal H}_{\rho,K_d}(x_0,\beta, \eta)$) is
included in the set of H\"older functions of order $\beta$ on $K_d$
(resp., around $x_0$).
\item Assuming $\be\le1$ is not restrictive since, for any $c>0$,
$\rho^c$ is also a quasi-metric.
\end{enumerate}
\end{rem}

The introduction of the logarithmic term
appears naturally when considering Gaussian random fields. Actually,
\cite{BEJARO} proves that for all $\beta\in(0,1]$, a large class of
elliptic Gaussian random fields $X_\beta$, including the famous
fractional Brownian fields, belongs a.s. to ${\mathcal H}_{\rho
,K_d}(x_0,\beta, 1/2)$ with $\rho$ the Euclidean distance
(see Theorem~1.3 in \cite{BEJARO}). Moreover, Xiao \cite{XiaoModC}
also gives some anisotropic examples of Gaussian fields
belonging a.s. to ${\mathcal H}_{\rho,K_d}(x_0,1, 1/2)$
for some anisotropic quasi-distance $\rho=\rho_{ {E}}$ associated
with $E$ a diagonal matrix (see Theorem~4.2 of \cite{XiaoAni}).
Finally, in \cite{BLS11}, we construct stable and Gaussian random
fields belonging a.s. to ${\mathcal H}_{\rho
_{x_0},K_d}(x_0,1-\varepsilon, 0)$ for some convenient $\rho_{x_0}$
(see Theorem~4.6 in \cite{BLS11}).

%s3 #&#
\section{Main results on conditionally sub-Gaussian series}
\label{MCGS}
%s3.1 #&#
\subsection{Local modulus of continuity}
\label{LMC}

In this section, we first give an upper bound of the local modulus of
continuity of $S$ defined by \eqref{SGSN} under the following local
assumption on the conditional parameter \eqref{VCS}.

%as2 #&#
%
\begin{hyp}
\label{Hx0} Let $x_0\in K_d$ with $K_d=\prod_{j=1}^d[a_j,b_j]\subset
\dR^d$.
Let us consider $\rho$ a quasi-metric on ${\mathbb R^d}$ satisfying equation
\eqref{controlrho}.
Assume that there exist an almost sure event $\Omega'$ and some random
variables
$\gamma>0$, $\beta\in(0,1]$, $\eta\in\dR$ and $C\in(0,+\infty)$
such that on $\Omega'$
\[
\forall x,y\in B(x_0,\gamma)\cap K_d, \qquad s\EM{{ (x,y
)}}\le C \rho(x,y)^{\beta}\bigl |\log\bigl(\rho (x,y)\bigr)\bigr |^{\eta},
\]
where we recall that the conditional parameter $s$ is given by \eqref{VCS}.
\end{hyp}

Note that the event $\Omega'$, the random variables $\gamma,\beta
,\eta$, $C$ and
the quasi-metric $\rho$ may depend on~$x_0$.

Let us now state the main result of this section on the modulus of
continuity. The main difference with \cite{Konot,BL09,BLS11} is that
we do not only consider the limit random field $S$ but obtain a uniform
upper bound in $N$ for the modulus of continuity of $S_N$.

%th3.1 #&#
%
\begin{theorem} \label{CVUS1}
Assume that Assumptions~\ref{HypTVg} and~\ref{Hx0} are fulfilled.
Then, almost surely, there exist $\gamma^*\in(0,\gamma)$ and $C \in
(0,+\infty)$ such that for all $x,y \in B(x_0,\gamma^*)\cap K_d$,
\[
\sup_{N\in\dN} {\EM{{\bigl\llvert S_N(x)-S_N(y)
\bigr\rrvert }}} \le C{\rho }(x,y)^{\beta} \EM{{\bigl\llvert \log\rho(x,y)
\bigr\rrvert }}^{\eta+1/2}.
\]
Moreover, almost surely $(S_N)_{N\in\dN}$ converges uniformly on
$B(x_0,\gamma^*)\cap K_d$ to $S$ and the limit $S$ belongs to
${\mathcal H}_{\rho,K_d}(x_0,\beta, \eta+1/2)$. In particular,
almost surely $S$ is continuous at $x_0$.
\end{theorem}

\begin{pf} See Appendix~\ref{LMC-A}.
\end{pf}

Strengthening Assumption~\ref{Hx0}, the uniform convergence and the
upper bound for the modulus of continuity are obtained on deterministic
set. Next corollary is obtained using some covering argument.

%co3.2 #&#
%
\begin{cor} Assume that Assumption~\ref{HypTVg} is fulfilled.
\label{ArgReco}
\begin{enumerate}[2.]
\item Assume that Assumption~\ref{Hx0} holds for any $x_0\in K_d$
with the same almost sure event $\Omega'$, the same random variables
$\beta$ and $\eta$, and the same quasi-metric $\rho$.
Then Theorem~\ref{CVUS1} holds replacing $B(x_0,\gamma^*)\cap K_d$ by
all the set $K_d$ and almost surely $S$ belongs to $ {\mathcal H}_{\rho
}(K_d,\beta, \eta+1/2)$.
\item Assume now that Assumption~\ref{Hx0} holds with a deterministic
$\gamma$. Then Theorem~\ref{CVUS1} holds
replacing $B(x_0,\gamma^*)\cap K_d$ by $B(x_0,\gamma)\cap K_d$
and almost surely $S$ belongs to $ {\mathcal H}_{\rho}(B(x_0,\gamma
)\cap K_d,\beta, \eta+1/2)$.
\end{enumerate}
\end{cor}
\begin{pf} See Appendix~\ref{LMC-A}.
\end{pf}

When considering $S$ an operator scaling Gaussian random field, note
that Li \textit{et al.} \cite{LiWangXiao13} proves that the upper
bound obtained by
Corollary~\ref{ArgReco} is optimal. Moreover for some Gaussian
anisotropic random fields, Xiao \cite{XiaoAni} also obtains a sample path
regularity in the stronger $L^p$-sense on whole the compact $K_d$. This
follows from an extension of the Garsia--Rodemich--Rumsey continuity
lemma Garsia \textit{et al.} \cite{GRR} or the minorization metric
method of Kwapie{\'n} and Rosi{\'n}ski \cite{KR04}.
This would be interesting to study if these results still hold when
considering a quasi-metric $\rho$ (and not a metric) and if they can
be applied to obtain the sample path regularity of $S$ in the stronger
$L^p$-sense on whole the compact $K_d$, strengthening the assumption on
the parameter $s$.

%s3.2 #&#
\subsection{Rate of almost sure uniform convergence}

This section is concerned with the rate of uniform convergence of the
series $\EM{{ (S_N )}}_{N\in\dN}$
defined by \eqref{SGSN}.
Under Assumption~\ref{HypTVg}, this series converges to $S$ and, for
any integer $N$, we
consider the
rest
\[
R_N(x)=S(x)-S_N(x)=\sum_{n=N+1}^{+\infty}W_n(x)g_n,
\qquad x\in K_d\subset \dR^d.
\]
Then, conditionally on $\EM{{ (W_n )}}_{n\ge1}$,
$R_N(x)-R_N(y)$ is a
sub-Gaussian random variable with parameter
%
%e9 #&#
%
\begin{equation}
\label{rN} r_N\EM{{ (x,y )}}=\EM{{ \Biggl(\sum
_{n=N+1}^{+\infty}\EM {{\bigl\llvert W_n(x)-W_n(y)
\bigr\rrvert }}^2 \Biggr)}}^{1/2},\qquad x, y\in
K_d.
\end{equation}
Observe that $R_0=S$ and that $r_0(x,y)=s(x,y)$.
To obtain a rate of uniform convergence for the sequence $\EM{{
(S_N )}}_{N\in\dN}$, the general assumption relies on a rate of
convergence for the sequence
$\EM{{ (r_N )}}_{N\in\dN}$.

%as3 #&#
%
\begin{hyp}
\label{Hx02} Let $x_0\in K_d$ with $K_d=\prod_{j=1}^d[a_j,b_j]\subset
\dR^d$ and let
$\rho$ be a quasi-metric on ${\mathbb R^d}$ satisfying \eqref{controlrho}.
Assume that there exist an almost sure event $\Omega'$, some random
variables $\gamma>0$, $\be\in(0,1]$, $\eta\in\dR$ and a positive
random sequence $(b(N))_{N\in\dN}$ such that
on $\Omega'$,
%
%e10 #&#
%
\begin{equation}
\label{VN} \forall N\in\dN, \forall x,y\in B(x_0,\gamma)\cap
K_d,\qquad r_N\EM {{ (x,y )}}\le b(N)
\rho(x,y)^{\beta}\bigl |\log\bigl(\rho (x,y)\bigr)\bigr |^{\eta}.
\end{equation}
\end{hyp}

Note that $\Omega'$, $\rho$ and the random variables $\gamma,\beta
,\eta$ and $b(N)$ may depend on $x_0$. Note also that since Assumption~\ref{Hx02} implies Assumption~\ref{Hx0}, according to Theorem~\ref
{CVUS1}, almost surely, there exists $\gamma^*\in(0,\gamma)$ such
that $R_N=S-S_N$ is continuous on $B(x_0,\gamma^*)$. The following
theorem precises the modulus of continuity of $R_N$ with respect to $N$
and a rate of uniform convergence.

%th3.3 #&#
%
\begin{theorem}
\label{CVUS2}
Assume that Assumptions~\ref{HypTVg} and~\ref{Hx02} are fulfilled.

\begin{enumerate}[2.]
\item Then, almost surely, there exists $\gamma^*\in(0,\gamma)$ and
$C \in(0,+\infty)$ such that for
\[
{\EM{{\bigl\llvert R_N(x)-R_N(y) \bigr\rrvert }}} \le
Cb(N) \sqrt{\log (N+2)} {\rho}(x,y)^{\beta
} \EM{{\bigl\llvert \log
\rho(x,y) \bigr\rrvert }}^{\eta+1/2}
\]
for all $N\in\dN$ and all $x,y \in B(x_0,\gamma^*)\cap K_d$.
\item Moreover, if almost surely, for all $N\in\dN$,
%
%e11 #&#
%
\begin{equation}
\label{Contx0} \EM{{\bigl\llvert R_N(x_0) \bigr\rrvert
}} \le b(N) \sqrt{\log(N+2)},
\end{equation}
then, almost surely, there exists $\gamma^*\in(0,\gamma)$ and $C\in
(0,+\infty)$ such that
\[
\EM{{\bigl\llvert R_N(x) \bigr\rrvert }}\le C{b(N) \sqrt{\log(N+2)}}
\]
for all $N\in\dN$ and all $x \in B(x_0,\gamma^*)\cap K_d$.
\end{enumerate}
\end{theorem}
\begin{pf} See Appendix~\ref{CVUS2-A}.
\end{pf}

An analogous of Corollary~\ref{ArgReco} holds for strengthening the
previous local theorem to get uniform results on $K_d$ or on $B(x_0,\ga
)\cap K_d$ when $\ga$ is deterministic.

%co3.4 #&#
%
\begin{cor}\label{ArgReco2} Assume that Assumptions~\ref{HypTVg} is fulfilled.
\begin{enumerate}[2.]
\item Assume that Assumption~\ref{Hx02} holds for any $x_0\in K_d$
with the same almost sure event $\Omega'$, the same random variables
$\beta$ and $\eta$, the same sequence $(b(N))_{N\in\dN}$ and the
same quasi-metric~$\rho$.
Then assertion 1 of Theorem~\ref{CVUS2} holds replacing $B(x_0,\gamma
^*)\cap K_d$ by all the set~$K_d$. If moreover, equation \eqref
{Contx0} is fulfilled for some $x_0\in K_d$, then
\[
\sup_{N\in\dN}\sup_{x\in K_d} \frac{\EM{{\llvert  R_N(x) \rrvert }}}{b(N) \sqrt
{\log(N+2)}}<+
\infty\qquad\mbox{almost surely.}
\]
\item Assume now that Assumption~\ref{Hx02} holds with a deterministic
$\gamma$. Then Theorem~\ref{CVUS2} holds
replacing $B(x_0,\gamma^*)\cap K_d$ by $B(x_0,\gamma)\cap K_d$.
\end{enumerate}
\end{cor}

%s4 #&#
\section{Shot noise series}
\label{SeSN}

%s4.1 #&#
\subsection{Preliminaries}

In this section, we consider the sequence of shot noise series defined by
\[
\forall N\in\dN, \forall\al\in K_1=[a,b]\subset(0,2),\qquad
S_N^*(\alpha)=\sum_{n=1}^NT_n^{-1/\al}X_n,
\]
where for all $n\ge1$, the random variable $T_n$ is the $n$th arrival
time of a Poisson process with intensity $1$ and $(X_n)_{n\ge1}$ is a
sequence of i.i.d. symmetric random variables, which is assumed
independent of $(T_n)_{n\ge1}$.
Let us first recall that $S_N^*(\alpha)$ converges almost surely to
$S^*(\al)$ an $\alpha$-stable random variable as soon as $X_n\in
L^\al$ (see, \cite{sam:taqqu}, e.g.). Under a strengthened
assumption on the integrability of $X_n$, rate of pointwise almost sure
convergence and rate of absolute convergence have also been given in
Theorems 2.1 and 2.2 of \cite{COLALE07}.

Since $\EM{{ (X_n )}}_{n\ge1}$ may not be a sequence of
sub-Gaussian
random variables, we cannot apply Section~\ref{MCGS} to the sequence
$\EM{{ (S_N^* )}}_{N\in\dN}$. However, due to symmetry of
$\EM{{ (X_n )}}_{n\ge1}$,
\[
\EM{{ (X_n )}}_{n\ge1} \stackrel{(d)} {=} \EM{{
(X_n g_n )}}_{n\ge1}
\]
with $\EM{{ (g_n )}}_{n\ge1}$ a Rademacher sequence
independent of $\EM{{ (X_n,T_n )}}_{n\ge1}$ and
Section~\ref{MCGS} allows to study
\[
S_N(\al)=\sum_{n=1}^N
W_n(\al) g_n \qquad\mbox{with $W_n(\al
):=T_n^{-1/\al} X_n$}.
\]
Moreover, in Theorems~\ref{CVUS1} and~\ref{CVUS2}, $S_N$
(resp., $R_N=S-S_N$) can be replaced by $S_N^*$ (resp.,
$R_N^*=S^*-S_N^*$), see the proof of next theorem for details. Then,
assuming that $X_n$ is sufficiently integrable, we obtain the uniform
convergence of $S_N^*$ on a deterministic compact interval
$K_1=[a,b]\subset(0,2)$ and a rate of uniform convergence. These
results, stated in the following theorem, strengthen Theorem~2.1 of
\cite{COLALE07} which deals with the pointwise rate of convergence.

%th4.1 #&#
%
\begin{theorem}\label{uCOLALE} For any integer $n\ge1$, let $T_n$ be
the $n$th arrival time of a Poisson process with intensity $1$. Let
$(X_n)_{n\ge1}$ be a sequence of i.i.d. symmetric random variables,
which is assumed independent of $(T_n)_{n\ge1}$.
Furthermore assume that
$
\mathbb{E}(|X_1|^{2p})<+\infty
$ for some $p>0$.
\begin{enumerate}[2.]
\item
Then, almost surely, for all $b\in(0,\min(2,2p))$ and for all $a\in
(0,b]$, the sequence of partial sums $\EM{{ (S_{N}^*
)}}_{N\in\dN}$
converges uniformly on $[a,b]$.
\item
Moreover, almost surely, for all $b\in(0,\min(2,2p))$ and for all
$a\in(0,b)$, for all $p'>0$ with
$1/p' \in(0,1/b-1/\min(2p,2))$,
\[
\sup_{N\in\dN}\sup_{\alpha\in[a,b]}N^{1/p'}\Biggl
\llvert \sum_{n=N+1}^{+\infty}T_n^{-1/\alpha}X_n
\Biggr\rrvert <+\infty.
\]
\end{enumerate}
\end{theorem}

\begin{pf} See Appendix~\ref{uCOLALE-A}.
\end{pf}
%
%s4.2 #&#
\subsection{Modulus of continuity and rate of convergence of shot
noise series}
\label{MCSN}

In this section, we focus on some shot noise series, which are
particular examples of conditionally sub-Gaussian series. For this
purpose, we assume that the following assumption is fulfilled.

%as4 #&#
%
\begin{hyp}
\label{HSN} Let $\EM{{ (T_n )}}_{n\ge1}$, $\EM{{
(V_n )}}_{n\ge1}$ and $\EM{{ (g_n )}}_{n\ge1}$ be
independent sequences satisfying the following conditions.
\begin{enumerate}[3.]
\item$\EM{{ (g_n )}}_{n\ge1}$ is a sequence of
independent complex-valued
symmetric sub-Gaussian random variables with parameter $s=1$.
\item$T_n$ is the $n$th arrival time of a Poisson process with
intensity $1$.
\item$\EM{{ (V_n )}}_{n\ge1}$ is a sequence of i.i.d.
complex-valued
random fields defined on $K_{d+1}\subset(0,2)\times\dR^d$.
\item For any $(\al,u)\in K_{d+1}$,
$
V_n(\al,u)\in L^{\al}$.
\end{enumerate}
\end{hyp}

For any integer $n\ge1$, we consider the complex-valued random field
$W_n$ defined by
%
%e12 #&#
%
\begin{equation}
\label{WnShotnoise} W_n\EM{{ (\al,u )}}:=T_n^{-1/\al}
V_n\EM{{ (\al ,u )}},\qquad(\al,u)\in K_{d+1}\subset(0,2)
\times\dR^d.
\end{equation}
Since $\EM{{\llvert  V_n(\alpha,u) \rrvert }}^{2}\in L^{\al/2}$ and
$\al/2\in(0,1)$,
according to Theorem~1.4.5 of \cite{sam:taqqu},
\[
\sum_{n=1}^{+\infty} \EM{{\bigl\llvert
W_n(\al,u) \bigr\rrvert }}^2=\sum
_{n=1}^{+\infty} T_n^{-2/\al} \EM{{
\bigl\llvert V_n(\alpha,u) \bigr\rrvert }}^{2}<+\infty
\qquad\mbox{almost surely}.
\]
Therefore, the independent sequences $\EM{{ (W_n )}}_{n\ge
1}$ and $\EM{{ (g_n )}}_{n\ge1}$ satisfy Assumption~\ref
{HypTVg}. Then, $S$ and
$(S_N)_{N\in\mathbb{N}}$ are well-defined on $K_{d+1}\subset
(0,2)\times\dR^d\subset\dR^{d+1}$ by \eqref{ChampS} and \eqref{SGSN}.
Before we study, the modulus of continuity of $S$ and the rate of
convergence of $(S_N)_{N\in\dN}$, let us state some remarks.

%re4.1 #&#
%
\begin{rem}
Assume that conditions 1--3 of Assumption~\ref{HSN} are fulfilled with
$\EM{{ (g_n )}}_{n\ge1}$
a sequence of i.i.d. random variables. Then
Remark~2.6 of \cite{Rosinskiu}
proves that condition 4 is a necessary and sufficient condition
for the almost sure convergence of $\EM{{ (S_N(\al,u)
)}}_{N\in\dN}$ for
each $(\al,u)\in K_{d+1}$. Note that by It\^o--Nisio theorem (see,
e.g., Theorem~6.1 of \cite{LedTal}), it is also a necessary and
sufficient condition for the convergence in distribution of the
sequence $\EM{{ (S_N\EM{{ (\al,u )}} )}}_{N\in
\dN}$. Then, condition 4
is not a strong assumption and is clearly essential to ensure that
$S\EM{{ (\al,u )}}$ is well-defined.
\end{rem}

%re4.2 #&#
%
\begin{rem}
\label{RepS} Assume that Assumption~\ref{HSN} is fulfilled with $\EM
{{ (g_n )}}_{n\ge1}$ a sequence of i.i.d. random variables.
Then, it is well
known that for each $\al\in(0,2)$, $S(\al,\cdot)$ is an $\al
$-stable symmetric random field, as field in variable $u$. In
Section~\ref{StableSE}, we will focus on $\al$-stable random fields defined
through a stochastic integral and see that, up to a multiplicative
constant, such a random field $X_\al$ has the same finite
distributions as $S(\al,\cdot)$ for a suitable choice of
$(g_n,V_n)_{n\ge1}$. The sample path regularity of $S$ in its variable
$\al$ is not needed to obtain an upper bound of the modulus of
continuity of $X_\al$. Nevertheless, this regularity is useful to deal
with multistable random fields (see Section~\ref{SeMStable}).
\end{rem}

The sequel of this section is devoted to simple criteria, based on some
moments of $V_n$, which ensure that Assumption~\ref{Hx02} (and then
Assumption~\ref{Hx0})
is fulfilled. More precisely, the results given below help us to give
simple conditions in order to get Assumption~\ref{Hx02} and \eqref
{Contx0} satisfied with
$b(N)=(N+1)^{-1/p'}$ for some convenient $p'>0$. Then,
all the results of Section~\ref{MCGS} hold.

%th4.2 #&#
%
\begin{theorem}\label{CrMom}
Assume that Assumption~\ref{HSN} is fulfilled with
$K_{d+1}=[a,b]\times\prod_{j=1}^d[a_j,b_j] \subset(0,2)\times\dR
^d$ and let $\rho$ be a quasi-metric on $\dR^{d+1}$ satisfying
equation \eqref{controlrho}. Assume also that for some $x_0\in
K_{d+1}$, there exist $r\in(0,+\infty)$, $\beta\in(0,1]$, $\eta\in
\dR$ and $p\in(b/2,+\infty)$ such that $\dE\EM{{ (\EM{{\llvert  V_1(x_0) \rrvert }}^{2p} )}}<+\infty$ and
%
%e13 #&#
%
\begin{equation}
\label{UBV} \dE\EM{{ \biggl(\EM{{ \biggl[\mathop{\sup_{x,y\in K_{d+1}}}_{0<\EM
{{\llVert  x-y\rrVert }}\le r}
\frac{\EM{{\llvert  V_1(x)-V_1(y)
\rrvert }}}{\rho(x,y)^{\beta}\EM{{\llvert  \log\rho(x,y) \rrvert }}^{\eta}} \biggr]}}^{2p} \biggr)}}<\infty.
\end{equation}
Let us recall that $S$ and $S_N$ are defined by \eqref{ChampS} and
\eqref{SGSN} with $W_n$ given by \eqref{WnShotnoise}.
\begin{enumerate}[2.]
\item Then, almost surely $\EM{{ (S_N )}}_{N\in\dN} $
converges uniformly
on $K_{d+1}$ and its limit $S$ belongs almost surely to $\mathcal
{H}_{\rho}\EM{{ (K_{d+1},\beta,\max\EM{{ (\eta,0
)}}+1/2 )}}$.
\item Moreover,
when $p'>0$ is such that $1/p'\in(0,1/b-1/\min(2p,2))$, almost surely
\[
\sup_{N\in\dN} N^{1/p'} \sup_{x\in K_{d+1}}
\EM{{\bigl\llvert S(x)-S_N(x) \bigr\rrvert }}<+\infty.
\]
\end{enumerate}
\end{theorem}

\begin{pf} See Appendix~\ref{SNMC-A}.
\end{pf}

%ex4.1 #&#
%
\begin{exple}
\label{expleG}
Assume that $V_1$ is a fractional Brownian field on $\dR^d$ with Hurst
parameter $H$.
Then \eqref{UBV} is satisfied for all $p>0$ with $\rho(x,y)=\|x-y\|$,
$\beta=H$ and $\eta=1/2$ (see, e.g., Theorem~1.3 of \cite{BEJARO}).
\end{exple}

Let us now
present a method (similar to those used in \cite{Kono,BL09,BLS11} to
bound some conditional variance) to establish \eqref{UBV}.

%pr4.3 #&#
%
\begin{prop}
\label{MG}
Let $x_0=\EM{{ (\al_0,u_0 )}}\in
K_{d+1}$ with $K_{d+1}=[a,b]\times\prod_{j=1}^d[a_j,b_j] \subset
(0,2)\times\dR^d$.
Let $V_1$ be a complex-valued random field defined on $K_{d+1}$. Assume
that there exists a random field $\EM{{ (\mathcal{G}(h)
)}}_{h\in
[0,+\infty)}$ with values in $[0,\infty)$ and such that
\begin{enumerate}[(iii)]
\item[(i)] there exists $\rho$ a quasi-metric on $\dR^{d+1}$
satisfying equation \eqref{controlrho} such that almost surely,
\[
\forall x,y\in K_{d+1},\qquad \EM{{\bigl\llvert V_1(x)-V_1(y)
\bigr\rrvert }}\le\mathcal{G}\bigl(\rho(x,y)\bigr);
\]
\item[(ii)] there exists $h_0\in(0,1]$ such that almost
surely, the function $h\mapsto\mathcal{G}(h)$ is monotonic on $[0,h_0]$;
\item[(iii)] there exist $p> b/2$ and some constants $\beta\in
(0,1]$, $\eta\in\dR$ and $C\in(0,\infty)$ such that for some
$\varepsilon>0$ and for $h>0$ small enough,
%
%e14 #&#
%
\begin{equation}
\label{IhTh} I(h):=\dE\EM{{ \bigl(\mathcal{G}(h)^{2p} \bigr)}}\le C
h^{2p\be
}\EM{{\llvert \log h \rrvert }}^{2p\EM{{ (\eta
-1/2p-\varepsilon )}}}.
\end{equation}
\end{enumerate}
Then, equation \eqref{UBV} holds for $r>0$ small enough.
\end{prop}

\begin{pf} See Appendix~\ref{SNMC-A}.
\end{pf}

%re4.3 #&#
%
\begin{rem}
If $ (V_n )_{n\ge1}$ is a sequence of independent symmetric
random variables, Theorem~\ref{CrMom} still holds replacing
$S_N\EM{{ (\al,u )}}$ (resp., $S\EM{{ (\al,u
)}}$) by
\[
S_N^*\EM{{ (\al,u )}}=\sum_{n=1}^{N}
T_n^{-1/\al} V_n\EM{{ (\al,u )}}\qquad \Biggl(
\mbox{resp., by } S^*\EM{{ (\al,u )}}=\sum_{n=1}^{+\infty}
T_n^{-1/\al} V_n\EM{{ (\al,u )}}\Biggr).
\]
In particular, following Example~\ref{expleG}, when Assumption~\ref
{HSN} is fulfilled with $V_1$ a
fractional Brownian field on ${\mathbb R^d}$ with Hurst parameter $H$,
assumptions of Theorem~\ref{CrMom} are fulfilled with $\rho_{d+1}$
the Euclidean distance on $\dR^{d+1}$, $\beta=H$ and $\eta=1/2$, on
any compact $(d+1)$-dimensional interval $K_{d+1}$. Especially, this
leads to an upper bound of the modulus of continuity of $S^*$
on any compact $(d+1)$-dimensional interval $K_{d+1}$.
Then for any fixed $\al_0\in(0,2)$, we also obtain
that the $\al_0$-stable random field $\EM{{ (S^*\EM{{ (\al
_0,u )}} )}}_{u\in
\dR^d}$ is in $\mathcal{H}_{\rho_d}\EM{{ (K_d,H,1 )}}$
for $\rho_d$ the
Euclidean distance on $\dR^d$ and for any compact set $K_d\subset\dR
^d$.
\end{rem}

%s4.3 #&#
\subsection{LePage random series representation}
\label{SeStable}

Representations in random series of infinitely divisible laws have been
studied in \cite{LePage1,LePage2}.
Such representations have been successfully used to study sample path
properties of some symmetric $\alpha_0$-stable random processes
($d=1$) and fields (see, e.g., \cite{Kono,BL09,BLS11,Dozzi11}).

Let us be more precise on the assumptions on the LePage series under study.

%as5 #&#
%
\begin{hyp} \label{HLep}
Let $\EM{{ (T_n )}}_{n\ge1}$ and $\EM{{ (g_n
)}}_{n\ge1}$ be as in Assumption~\ref{HSN}. Let ${\EM{{ (\xi_n )}}}_{n \ge1}$
be a sequence
of i.i.d. random
variables with common law
\[
\mu\EM{{ (\mathrm{d}\xi )}}=m\EM{{ (\xi )}} \nu\EM {{ (\mathrm{d}\xi )}}
\]
equivalent to a $\sigma$-finite measure $\nu$ on $(\dR^d,{\mathcal
B}(\dR^d))$ (that is such that $m(\xi)>0$ for $\nu$-almost every
$\xi$).
This sequence is independent from $\EM{{ (g_n,T_n )}}_{n\ge1}$.
Moreover, we consider
\[
V_n\EM{{ (\al,u )}}:=f_{\al}\EM{{ (u,\xi_n
)}} m\EM{{ (\xi_n )}}^{-1/\al},
\]
where for any $\al\in K_1\subset(0,2)$, $f_{\al} \dvtx K_d\times
\dR
^d\rightarrow\dC$ is a deterministic function such that
\[
\forall u\in K_d\subset\dR^d,\qquad \int
_{\dR^d} \EM{{\bigl\llvert f_{\al}\EM{{ (u,\xi )}}
\bigr\rrvert }}^\al\nu\EM{{ (\mathrm{d}\xi )}}<+\infty.
\]
\end{hyp}

Under this assumption,
Assumption~\ref{HSN} is fulfilled with $K_{d+1}=K_1\times K_d$. Then,
emphasizing the dependence on the function $m$, $S_{m,N}$ and $S_m$ are
well defined on $K_{d+1}$ by \eqref{SGSN} and \eqref{ChampS} with
$W_n$ given by \eqref{WnShotnoise}. In particular,
%
%e15 #&#
%
\begin{eqnarray}
\label{SLePage}  &&S_m\EM{{ (\al,u )}}=\sum
_{n=1}^{+\infty} T_n^{-1/\al}
f_{\al}\EM{{ (u,\xi_n )}} m\EM{{ (\xi_n
)}}^{-1/\al} g_n,
\nonumber
\\[-8pt]
\\[-8pt]
&&\quad (\alpha,u)\in K_{d+1}:=K_1
\times K_d\subset(0,2)\times\dR^d.\nonumber
\end{eqnarray}
Under appropriate assumptions on $f_{\al}$ and $m$, the previous
sections state the uniform convergence of the series, give a rate of
convergence and some results on regularity for
$S_m$.
Precise results on regularity of $S_m$ may be obtained using the
following proposition, which states that the finite distributions of
$S_m$ does not depend on the choice of the $\nu$-density $m$.

%pr4.4 #&#
%
\begin{prop}\label{eq:Lepage} Assume that Assumption~\ref{HLep} is
fulfilled and let $S_m$ be defined by \eqref{SLePage}.
Let $(\tilde{\xi}_n)_{n\ge1}$ be a sequence of i.i.d. random
variables with common law
$
\tilde{\mu}\EM{{ (\mathrm{d}\xi )}}=\tilde{m}\EM{{ (\xi
 )}} \nu\EM{{ (\mathrm{d}\xi )}}
$
equivalent to $\nu$. Assume that the sequences $(\tilde{\xi
}_n)_{n\ge1}$, $\EM{{ (g_n )}}_{n\ge1}$ and $\EM{{
(T_n )}}_{n\ge1}$ are independent.
\begin{enumerate}[2.]
\item Then, $S_m\stackrel{\mathit{fdd}}{=}S_{\tilde{m}}$, where $\stackrel
{\mathit{fdd}}{=} $ means equality of finite distributions. In other words,
\[
\EM{{ \bigl(S_m(\al,u) \bigr)}}_{(\al,u)\in K_{d+1}}\stackrel {\mathit{fdd}} {=}
\EM{{ \Biggl(\sum_{n=1}^{+\infty}
T_n^{-1/\al} f_{\al
}\EM{{ (u,\tilde{
\xi}_n )}} \tilde{m}\EM{{ (\tilde {\xi}_n
)}}^{-1/\al} g_n \Biggr)}}_{(\al,u)\in K_{d+1}}.
\]
\item Assume moreover that for $\nu$-almost every $\xi\in\dR^d$,
the map $\EM{{ (\al,u )}}\mapsto f_\al\EM{{ (u,\xi
 )}}$ is continuous on the
compact set $K_{d+1}\subset(0,2)\times\dR^d$. Let us consider $\rho
$ a quasi-metric on $\dR^{d+1}$, $\beta\in(0,1]$ and $\eta\in\dR
$. Then, $S_m$ belongs almost surely in $\mathcal{H}_{\rho
}(K_{d+1},\beta,\eta)$ if and only if $S_{\tilde{m}}$ does.
\end{enumerate}
\end{prop}

\begin{pf}
See Appendix~\ref{eq:Lepage-A}.
\end{pf}

In particular,
when studying the sample path properties of $S_m$, this result allows
us to replace $m$ by an other function $\tilde{m}$
so that the regularity of $S_m$ may be deduced from the regularity of
$S_{\tilde{m}}$. For example, replacing $m$ by $m_{x_0}$
depending on $x_0$ this may lead to a more precise bound for the
modulus of continuity of $S_m$ around $x_0$
(see, e.g., Example~\ref{expleMO}).\vadjust{\goodbreak}

%s5 #&#
\section{Applications}\label{Applications}

%s5.1 #&#
\subsection{\texorpdfstring{$\alpha$}{alpha}-stable isotropic random fields}
\label{StableSE}

Let us fix $\alpha=\alpha_0\in(0,2)$ and assume that Assumption~\ref
{HLep} is fulfilled with $g_n$ some isotropic complex random variables. Then,
the proof of Proposition~\ref{eq:Lepage} (see Section~\ref
{eq:Lepage-A}) allows to compute
the characteristic function of the isotropic $\al_0$-stable random
field $S_m\EM{{ (\al_{0},\cdot )}}=\EM{{ (S_m\EM
{{ (\al_0,u )}} )}}_{u\in K_d}$,
which leads to
\[
S_m\EM{{ (\al_0,\cdot )}}\stackrel{\mathrm{fdd}} {=}
d_{\al_0} \EM{{ \biggl(\int_{\dR ^d} f_{{\al_0}}
\EM{{ (u,\xi )}} M_{\al_0}\EM{{ (\mathrm{d}\xi )}} \biggr)}}_{u\in K_d},
\]
with $M_{\al_0}$ a complex isotropic $\al_0$-stable random measure on
$\dR^d$ with control measure $\nu$ and
%
%e16 #&#
%
\begin{equation}
\label{defcste} d_{\al_0}=\dE\EM{{ \bigl(\EM{{\bigl\llvert \Re\EM{{
(g_1 )}} \bigr\rrvert }}^{\al_0} \bigr)}}^{1/{\al_0}}\EM{{
\biggl(\frac
{1}{2\pi}\int_{0}^{2\pi}\EM{{\bigl
\llvert \cos\EM{{ (\theta )}} \bigr\rrvert }}^{\al_0} \,\mathrm{d}\theta
\biggr)}}^{-1/{\al_0}}\EM {{ \biggl(\int_{0}^{+\infty}
\frac{\sin\EM{{ (\theta
)}}}{\theta^{\al_0}}\,\mathrm{d}\theta \biggr)}}^{1/{\al_0}}.
\end{equation}
When $\nu$ is a finite measure (resp., the Lebesgue measure),
this stochastic integral representation of $S_m\EM{{ (\al
_0,\cdot )}}$
has been provided in \cite{sam:taqqu,Marcus} (resp., \cite{Kono,BL09}).

Let us note that assumptions of Theorem~\ref{CrMom} and Proposition~\ref{MG} can be stated in term of the deterministic kernel $f_{{\al
_0}}$ to obtain an upper bound of the modulus of continuity of $S_m$.
In general, well-choosing $m_{u_0}$ and applying Proposition~\ref
{eq:Lepage}, we obtain a more precise upper bound of the modulus of
continuity of $S_{m}\EM{{ (\al_0,\cdot )}}$ around $u_0$,
which also
holds for a modification of the random field
%
%e17 #&#
%
\begin{equation}
\label{RepInt} X_{\al_0} =\EM{{ \biggl(\int_{\dR^d}
f_{{\al_0}}\EM{{ (u,\xi )}} M_{\al _0}\EM{{ (\mathrm{d}\xi )}}
\biggr)}}_{u\in K_d}.
\end{equation}
To illustrate how the previous sections can be applied to study the
field $X_{\al_0}$, which is defined through a stochastic integral and
not a series, let us focus on the case of harmonizable stable random
fields. More precisely, we consider
%
%e18 #&#
%
\begin{equation}
\label{NoyHarmo} f_{{\al_0}}\EM{{ (u,\xi )}} =\EM{{ \bigl(\mathrm
{e}^{\mathrm{i}\langle u,\xi\rangle} -1 \bigr)}} \psi_{\al_0}\EM {{ (\xi )}},\qquad
\forall(u,\xi)\in\dR^d\times\dR^d,
\end{equation}
with
$\psi_{\al_0}\dvtx\dR^d\rightarrow\dC$ a Borelian function such that
\[
\int_{\dR^d}\min\bigl(1,\|\xi\|^{\al_0}\bigr) \EM{{
\bigl\llvert \psi_{\al
_0}\EM{{ (\xi )}} \bigr\rrvert }}^{\al_0}
\nu(\mathrm{d}\xi) <+\infty.
\]
Note that,
since this assumption does not depend on $u$, the random
field $X_ {\alpha_0}$ may be defined on the whole space $\dR^d$.
For the sake of simplicity, in the sequel, we consider the case where
$\nu$ is the Lebesgue measure and first focus on a random field
$X_{\al_0}$ which behaves as operator scaling random fields studied in
\cite{OSSRF}.

%pr5.1 #&#
%
\begin{prop}
\label{SIOS} Let $\al_0\in(0,2)$ and let $X_{\al_0}$ be defined by
\eqref{RepInt} with $\nu$ the Lebesgue measure on $\dR^d$.
Let $E$ be a real matrix of size $d\times d$ whose eigenvalues have
positive real parts. Let $\tau_{ {E}}$ and $\tau_{ {{E^t}}}$ be
functions as introduced in Example~\ref{tauE} and let us set
$
q(E)=\operatorname{trace}(E)$ and
$a_1=\min_{\la\in\operatorname{Sp}(E)}\Re\EM{{ (\la )}}
$
with $ \operatorname{Sp}(E)$ the spectrum of $E$, that is, the set of the
eigenvalues of $E$. Assume that there exist some finite positive
constants $c_\psi$, $A$ and $\beta\in(0,a_1)$ such that
%
%e19 #&#
%
\begin{equation}
\label{CondHS} \EM{{\bigl\llvert \psi_{\al_0}\EM{{ (\xi )}} \bigr\rrvert
}} \le c_\psi\tau_{ {{E^t}}}(\xi )^{-\beta-q(E)/\al_0}, \qquad
\mbox{for almost every } \|\xi\|>A.
\end{equation}
Then, there exists a modification $X_{\al_0}^*$ of $X_{\al_0}$ such
that almost surely, for any $\varepsilon>0$, for any non-empty compact
set $K_d\subset\dR^d$,
\[
\mathop{\sup_{u,v\in K_d}}_{u\ne v}\frac{\EM{{\llvert  X_{\alpha
_0}^*(u)-X_{\alpha_0}^*(v) \rrvert }}}{\tau_{ {E}}\EM{{
(u-v )}}^{\beta} \EM{{ [\log\EM{{ (1+\tau_{
{E}}\EM{{ (u-v )}}^{-1} )}} ]}}^{\varepsilon
+1/2+1/\alpha
_0}}<+\infty.
\]
\end{prop}

%re5.1 #&#
%
\begin{rem}
\label{remComp} The quasi-metric $\EM{{ (x,y )}}\mapsto
\tau_{ {E}}\EM{{ (x-y )}}$ may not fulfill equation
\eqref{controlrho} since the
eigenvalues of $E$ may not be greater than 1. Nevertheless, the
quasi-metric $ \EM{{ (x,y )}}\mapsto\tau_{ {{E/a_1}}}\EM
{{ (x-y )}}$ does
and the conclusion with $\tau_{{E}}$ in the previous proposition
then follows from the comparison
\[
\forall\xi\in\dR^d,\qquad c_{1}\tau_{ {E}}(
\xi)^{a_1}\le\tau_{{E/a_1}}(\xi)\le c_{{2}}
\tau_{ {E}}(\xi)^{a_1}
\]
with $c_1,c_2$ two finite positive constants.
\end{rem}
\begin{pf} See Appendix~\ref{SIOS-A}.
\end{pf}

An upper bound for the modulus of continuity of such harmonizable
random fields is also obtained in
\cite{XiaoModC}. This upper bound is given in term of the Euclidean
norm and then does not take into account the anisotropic behavior of
$X_{\al_0}$. Even when $\tau_{ {E}}$ is the Euclidean norm, our
result is a little more precise than the one of \cite{XiaoModC}. The
difference is only in the power of the logarithmic term.

Let us now give some examples. We keep the notation of the previous
proposition and the eigenvalues of the matrix $E$ have always positive
real parts.

%ex5.1 #&#
%
\begin{exple}[(Operator scaling random fields \cite{OSSRF})] Let
$\psi\dvtx\dR^d\rightarrow[0,\infty)$ be an $E^t$-homogeneous
function, which means that
\[
\forall c\in(0,+\infty), \forall\xi\in\dR^d,\qquad \psi\EM {{
\bigl(c^{E^t} \xi \bigr)}}=c\psi\EM{{ (\xi )}},
\]
where $c^{E^t}=\exp\EM{{ (E^t \log c )}}$. Let us assume
that $\psi$ is a
continuous function such that $\psi(\xi)\ne0$ for $\xi\ne0$. Then
we consider the function $\psi_{\al_0}\dvtx\dR^d\rightarrow
[0,+\infty
]$ defined by
\[
\psi_{\al_0}\EM{{ (\xi )}}=\psi(\xi)^{-H-q(E)/\al_0}.
\]
The random field $X_{\al_0}$, associated with $\psi_{\al_0}$ by
\eqref{RepInt} and \eqref{NoyHarmo}, is well-defined and is
stochastically continuous if and only if $H\in(0,a_1)$. Then, let us
now fix $H\in(0,a_1)$.
Since $\psi_{\al_0}$ is $E^t$-homogeneous, one easily checks that
there exists $c_\psi\in(0,+\infty)$ such that
\[
\forall\xi\in\dR^d,\qquad\psi_{\al_0}\EM{{ (\xi )}}\le
c_\psi\tau_{
{{E^t}}}(\xi)^{-H-q(E)/\al_0} .
\]
Then, the assumptions of Proposition~\ref{SIOS} are fulfilled with
$\be=H$.
The corresponding conclusion was stated in Theorem~5.1 of \cite{BL09}
when $H=1$ and $ a_1>1$, which is enough to cover the general case
using Remark~2.1 of \cite{BL09}.
\end{exple}

%ex5.2 #&#
%
\begin{exple}[(Anisotropic Riesz--Bessel $\boldsymbol{\al}$-stable
random fields)]
\label{expleBessel} Let us consider
\[
\psi_{\al_0}\EM{{ (\xi )}}=\frac{1}{\tau_{
{{E^t}}}(\xi)^{2\be
_1/\al_0}\EM{{ (1+\tau_{ {{E^t}}}(\xi)^{2} )}}^{\beta
_2/\al_0}},\qquad\xi \in
\dR^d\backslash\EM{{ \{0 \}}}
\]
with two real numbers $\be_1$ and $\be_2$.
Assuming that
\[
\frac{q(E)}{2}< \be_1+\be_2 \quad\mbox{and}\quad
\be_1< \frac
{q(E)}{2}+\frac{\al_0 a_1}{2},
\]
the random field $X_{\al_0}$ is well-defined by \eqref{RepInt}.
When $\tau_{ {{E^t}}}$ is the Euclidean norm, this random field has
been introduced in \cite{XiaoModC}
to generalize the Gaussian fractional Riesz--Bessel motion \cite{Anh}.

We distinguish two cases. If
$
\beta_1+\beta_2<\frac{q(E)}{2}+\frac{\al_0 a_1}{2}$, Proposition~\ref{SIOS} can be applied with $\beta=\frac{2(\beta
_1+\beta_2)-q(E)}{\al_0}$. Otherwise, Proposition~\ref{SIOS} can be
applied for any $\beta\in(0,a_1)$.
\end{exple}

Random fields defined by \eqref{RepInt} have stationary increments so
that their regularity on $K_d$ does not depend on the compact set
$K_d$. To avoid this feature, one can consider non-stationary
generalizations by substituting $\psi_{\alpha_0}$ by a function that
also depends on $u\in K_d$.
More precisely, we can consider
%
%e20 #&#
%
\begin{equation}
\label{RepInt2} X_ {\alpha_0}{=} \EM{{ \biggl(\int_{\dR^d}
\EM{{ \bigl(\mathrm {e}^{\mathrm{i}\langle u,\xi\rangle} -1 \bigr)}} \psi_{\al_0}\EM{{ (u,
\xi )}} {M}_{\al_0}\EM{{ (\mathrm{d}\xi )}} \biggr)}}_{
u\in K_d}
\end{equation}
with ${M}_{\al_0}$ a complex isotropic $\al_0$-stable random measure
with Lebesgue control measure and $\psi_{\al_0}$ a Borelian function
such that, for all $u\in K_d$,
\[
\int_{\dR^d}\EM{{\bigl\llvert \mathrm{e}^{\mathrm{i}\langle u,\xi\rangle} -1
\bigr\rrvert }}^{\al_0} \EM{{\bigl\llvert \psi_{\al_0}\EM{{ (u,\xi
)}} \bigr\rrvert }}^{\al_0}\,\mathrm{d}\xi<+\infty.
\]
Under some conditions on $\psi_{\al_0}$, when considering the local
behavior of $X_{\alpha_0}$ around a point $u_0$ one can
conveniently choose a Lebesgue density $m_{u_0}$ to obtain an upper
bound of the modulus of continuity of the shot noise series
$S_{m_{u_0}}\EM{{ (\al_0,\cdot )}}$ given by \eqref
{SLePage} with
\[
f_{{\al_0}}\EM{{ (u,\xi )}}=\EM{{ \bigl(\mathrm {e}^{\mathrm{i}\langle u,\xi\rangle} -1
\bigr)}} \psi_{\al_0}\EM{{ (u,\xi )}}.
\]
For the sake of conciseness, let us illustrate this with multi-operator
random fields, which have already been studied in \cite{BLS11}.

%ex5.3 #&#
%
\begin{exple}[(Multi-operator scaling $\boldsymbol{\al}$-stable
random fields)]\label{expleMO}
In \cite{BLS11}, we consider
$E$ a function defined on $\dR^d$ with values in the set of real
matrix of size $d\times d$ whose eigenvalues have real parts greater
than $1$ and $\psi\dvtx{\mathbb R^d}\times{\mathbb R^d}\rightarrow
[0,+\infty)$ a
continuous function such that for any $u\in{\mathbb R^d}$,
$\psi(u,\cdot)$ is homogeneous with respect to $E(u)^t$, that is,
\[
\psi\bigl(u,c^{E(u)^t}\xi\bigr)=c\psi(u, \xi),\qquad\forall\xi\in{\mathbb
R^d}, \forall c>0.
\]
Under convenient regularity assumptions on $\psi$ and $E$,
the ${\al_0}$-stable random field $X_{\al_0}$ is well-defined by
\eqref{RepInt2} setting
\[
\psi_{\al_0}(u,\xi)=\psi(u, \xi)^{-1-q(E(u))/\al_0} \qquad\mbox{with } q
\bigl(E(u)\bigr)=\operatorname{trace}\bigl(E(u)\bigr).
\]
Let $K_d=\prod_{j=1}^d[a_j,b_j]\subset\dR^d$ and $u_0\in K_d$. Let
us set $K_{d+1}=\{\al_0\}\times K_d$ and consider the quasi-metric
$\rho$ defined on $\dR^{d+1}$ by
\[
\rho\bigl((\alpha,u),\bigl(\alpha',v\bigr)\bigr)=\bigl |\alpha-
\alpha'\bigr |+\tau_{ {{E(u_0)}}}(u-v)
\]
for all $(\alpha,u),(\alpha',v)\in\dR\times{\mathbb R^d}$, which clearly
satisfies equation \eqref{controlrho}. Then, under assumptions of
\cite{BLS11}, there exists a Lebesgue density $m_{u_0}>0$ a.e. such
that Assumption~\ref{Hx0} holds for $S_{m_{u_0}}$ on $K_{d+1}$ with $\eta=0$ and all
$\beta\in(0,1)$, adapting similar arguments as in Proposition~\ref
{SIOS} (see Lemma~4.7 of \cite{BLS11}).
Therefore, following a part of the proof of Proposition~\ref{SIOS},
there exists a modification $X_{\al_0}^*$ of $X_{\al_0}$ such that
almost surely,
\[
\lim_{r\downarrow0} \mathop{\sup_{u,v\in B(u_0,r)\cap K_d}}_{u\ne
v}
\frac{\EM{{\llvert  X_{\al _0}^*(u)-X_{\al_0}^*(v) \rrvert }}}{\tau_{{E(u_0)}}(u-v)^{1-\veps}}<+\infty
\]
for any $\veps\in(0,1)$. This is Theorem~4.6 of \cite{BLS11}.
\end{exple}

For the sake of conciseness, we do not develop other examples.
Nevertheless, let us mention that our results can also be applied to
harmonizable fractional $\al$-stable sheets or even to operator stable
sheets. In particular, this improves the result stated in \cite{Xiao}
for fractional $\al$-stable sheets. Note that we can also deal with
real symmetric measure $W_\al$.

%s5.2 #&#
\subsection{Multistable random fields}
\label{SeMStable}

Multistable random fields have first been introduced in \cite{FLV09}
and then studied in \cite{FLGLV09}. Each marginal $X(u)$ of such a
random field is a stable random variable but its stability index is
allowed to depend on the position $u$.

Generalizing the class of multistable random fields introduced in \cite
{LGLV12}, we consider a multistable random field defined by a LePage
series. More precisely, under Assumption~\ref{HLep}, we consider
%
%e21 #&#
%
\begin{equation}
\label{multi} \tilde{S}_m\EM{{ (u )}}=\sum
_{n=1}^{+\infty} T_n^{-1/\al
(u)}
f_{{\al
\EM{{ (u )}}}}\EM{{ (u,\xi_n )}} m\EM{{ (\xi_n
)}}^{-1/\al(u)} g_n,\qquad u \in K_d,
\end{equation}
where $\al\dvtx K_d\rightarrow(0,2)$ is a function.
Then since
$
\tilde{S}_m\EM{{ (u )}}=S_m\EM{{ (\al(u),u )}}
$
with $S_m$ defined by \eqref{SLePage}, we deduce from Section~\ref{SeSN}
an upper bound for the modulus of continuity of $\tilde{S}$.
In particular, assuming that $\al$ is smooth enough, we obtain the
following theorem.

%pr5.2 #&#
%
\begin{prop}\label{ModCMS}
Let $K_d=\prod_{j=1}^{d}[a_j,b_j]\subset\dR^d$. Let us choose
$u_0\in K_d$. Let $\tilde{\rho}$ be a quasi-metric on $\dR^d$
satisfying equation \eqref{controlrho} and let $\al\dvtx
K_d\rightarrow
(0,2)$ belongs to $\mathcal{H}_{{\tilde{\rho}}}\EM{{
(K_d,1,0 )}}$. Let
us set
\[
a=\min_{{K_{ d}}}\al,\qquad b=\max_{{K_{ d}}} \al
\quad\mbox {and}\quad K_{1}=[a,b] \subset(0,2)
\]
and consider the quasi-metric $\rho$ defined on $\dR\times\dR^d$ by
\[
\rho\EM{{ \bigl(\EM{{ (\al,u )}},\EM{{ \bigl(\al ',v \bigr)}}
\bigr)}}=\EM{{\bigl\llvert \al-\al' \bigr\rrvert }}+ \tilde {\rho}
\EM{{ (u,v )}}.
\]
Assume that Assumption~\ref{HLep} is fulfilled and that equation
\eqref{UBV} holds on $K_{d+1}=[a,b]\times K_d$ for some $p>b/2$,
$\beta\in(0,1]$ and $\eta\in\dR$. Assume also that
\[
\dE\EM{{ \bigl(\EM{{\bigl\llvert V_1\bigl(\al(u_0),u_0
\bigr) \bigr\rrvert }}^{2p} \bigr)}}=\int_{{\mathbb R^d}}\EM{{
\bigl\llvert f_{{\al \EM{{ (u_0
)}}}}\EM{{ (u_0,\xi )}} \bigr\rrvert
}}^{2p} m\EM{{ (\xi )}}^{1-2p/\al(u_0)}\,\mathrm{d}\xi <+\infty.
\]
Let ${S}_{m,N}$ be defined by \eqref{SGSN} with $W_n\EM{{ (\al
,u )}}=T_n^{-1/\al} f_{{\al}}\EM{{ (u,\xi_n )}}
m\EM{{ (\xi_n )}}^{-1/\al}$
and let $\tilde{S}_{m,N}(u)=S_{m,N}(\al(u),u)$.
\begin{enumerate}[2.]
\item Then, almost surely, $(\tilde{S}_{m,N})_{N\in\dN}$ converges
uniformly on $K_d$ to $\tilde{S}_m$ and {almost surely} the limit
$\tilde{S}_m$ belongs to $ {\mathcal H}_{\tilde{\rho}}(K_d,\beta
,\max\EM{{ (\eta,0 )}}+1/2)$.
\item Moreover,
for all $p'>0$ such that $1/p'\in(0,1/b-1/\min(2p,2))$,
\[
\sup_{N\in\dN} N^{1/p'}\sup_{u\in K_d}\EM{{
\bigl\llvert \tilde{S}_m(u)-\tilde {S}_{m,N}(u) \bigr
\rrvert }}<+\infty.
\]
\end{enumerate}
\end{prop}

\begin{pf} See Appendix~\ref{SeMStable-A}.
\end{pf}

%re5.2 #&#
%
\begin{rem} Let us recall that $\tilde{S}_{m}\in{\mathcal H}_{\tilde
{\rho}}(K_d,\beta,\max\EM{{ (\eta,0 )}}+1/2)$ if and
only if $\tilde
{S}_{\tilde{m}}\in{\mathcal H}_{\tilde{\rho}}(K_d,\beta,\max\EM
{{ (\eta,0 )}}+1/2)$, with $\tilde{m}$ an other $\nu
$-density equivalent
to $\nu$, by Proposition~\ref{eq:Lepage}.
\end{rem}

To illustrate the previous proposition, we only focus on multistable
random fields obtained replacing in a LePage series representation of
an harmonizable operator scaling stable random field the index $\al$
by a function. Many other examples can be given, such as multistable
anisotropic Riesz--Bessel random fields or the class of linear
multistable random fields defined in \cite{FLGLV09}.

%co5.3 #&#
%
\begin{cor}[(Multistable versions of harmonizable operator
scaling random fields)]\label{MSV}
Let $E$ be a real matrix of size $d\times d$ such that $\min_{\la\in
\operatorname{Sp} E}\Re\EM{{ (\la )}}>1$. Let us
consider ${{\rho_{ {E}}}}$
and ${{\tau_{ {E}}}}$ as defined in Example~\ref{tauE}. Let us also
consider $\psi\dvtx\dR^d\rightarrow[0,\infty)$ a continuous,
$E^t$-homogeneous function such that $\psi(\xi)\ne0$ for $\xi\ne
0$. Then we set
\[
f_{\al}\EM{{ (u,\xi )}}=\EM{{ \bigl(\mathrm {e}^{\mathrm{i}\langle u, \xi\rangle} -1
\bigr)}}\psi(\xi)^{-1-q(E)/\al}
\]
with $q(E)=\operatorname{trace}(E)$. Let $m$ be a Lebesgue density a.e.
positive on ${\mathbb R^d}$,
$\EM{{ (\xi_n,T_n,g_n )}}_{n\ge1}$ be as in Assumption~\ref{HLep} with
$\nu$ the Lebesgue measure and consider a function $\al\dvtx\dR
^d\rightarrow(0,2)$. Therefore, the multistable random field $\tilde
{S}_m$ is well-defined by \eqref{multi} on the whole space ${\mathbb R^d}$.
Moreover if $\al\in\mathcal{H}_{{\rho_{ {E}}}}\EM{{
({\mathbb R^d},1,0 )}}$,
then for any $u_0\in\dR^d$ and $\veps>0$, there exists $r\in(0,1]$
such that
almost surely
\[
\mathop{\sup_{u,v\in B(u_0,r)}}_{u\ne v}\frac{\EM{{\llvert  \tilde
{S}_m\EM{{ (u )}}-\tilde{S}_m(v) \rrvert }}}{\tau
_{{E}}\EM{{ (u-v )}}\EM{{\llvert  \log\tau_{ {E}}\EM
{{ (u-v )}} \rrvert }}^{1/\al(u_0)+1/2+\veps}}<+
\infty.
\]
\end{cor}

\begin{pf} See Appendix~\ref{SeMStable-A}.
\end{pf}

%re5.3 #&#
%
\begin{rem}
In particular, when $E=\operatorname{Id}$, $\tau_{{E}}$ is the Euclidean
norm and we obtain an upper bound of the modulus of continuity of
multistable versions of fractional harmonizable stable fields.
\end{rem}

\begin{appendix}
%s6 #&#
\section{Proof of Proposition \texorpdfstring{\protect\ref{SGP}}{2.1}}
\label{Prel-A}

The proof of Proposition~\ref{SGP} is based on the following lemma.

%le6.1 #&#
%
\begin{lemma}
\label{DecSG}
If $Z$ is a complex-valued sub-Gaussian random variable with parameter
$s\in(0,+\infty)$, then for all
$t\in(0,+\infty)$, $ \dP\EM{{ (\EM{{\llvert  Z \rrvert }}>
t )}} \le4\mathrm{e}^{-\frace
{t^2}{8s^2}}$.
\end{lemma}
\begin{pf} Let $t\in(0,+\infty)$.
Since $Z$ is sub-Gaussian with parameter $s$, $\Re\EM{{ (Z
)}}$ and $\Im
\EM{{ (Z )}}$ are real-valued sub-Gaussian random variables
with parameter
$s$. Then applying Proposition~4 of \cite{Kahane60},
\[
\dP\EM{{ \bigl(\EM{{\llvert Z \rrvert }}> t \bigr)}}\le\dP\EM {{ \biggl(\EM{{
\bigl\llvert \Re\EM{{ (Z )}} \bigr\rrvert }}> \frac
{t}{2} \biggr)}}+\dP
\EM{{ \biggl(\EM{{\bigl\llvert \Im\EM{{ (Z )}} \bigr\rrvert }}> \frac{t}{2}
\biggr)}} \le4 \exp\EM{{ \biggl(-\frac{t^2}{8 s^2} \biggr)}},
\]
which
concludes the proof.
\end{pf}

\renewcommand{\theequation}{\arabic{equation}}
\setcounter{equation}{21}

Let us now prove Proposition~\ref{SGP}.
\begin{pf*}{Proof of Proposition~\ref{SGP}}
Let $t\in(0,+\infty)$. Since Proposition~\ref{SGP} is straightforward
if $a=0$, we assume that $a\ne0$.
Since the sequence $\EM{{ (g_n )}}_{n\ge1}$ is symmetric,
by the L\'evy
inequalities (see Proposition~2.3 in \cite{LedTal}), for any $M\in\dN
\backslash\EM{{ \{0 \}}}$,
\[
\dP\EM{{ \Biggl(\sup_{1\le P\le M} \EM{{\Biggl\llvert \sum
_{n=1}^P a_n g_n
\Biggr\rrvert }}> t\EM{{\llVert a\rrVert }}_{{\ell^2}} \Biggr)}}\le 2\dP\EM{{
\Biggl(\EM{{\Biggl\llvert \sum_{n=1}^M
a_n g_n \Biggr\rrvert }}> t\EM {{\llVert a\rrVert
}}_{{\ell^2}} \Biggr)}}.
\]
We now prove that $ {\sum_{n=1}^M a_n g_n}$ is sub-Gaussian. By
independence of the random variables $g_n$ and since each $g_n$ is
sub-Gaussian with parameter $s=1$,
%
%e22 #&#
%
\begin{equation}
\label{SGSF} \forall z\in\dC,\qquad\dE\EM{{ \bigl(\mathrm{e}^{\Re\EM{{
(\overline{z} \sum _{n=1}^M a_n g_n )}}}
\bigr)}}=\prod_{n=1}^M\dE\EM{{ \bigl(
\mathrm{e}^{\Re\EM{{ (\overline{z}
a_n g_n )}}} \bigr)}} \le\prod_{n=1}^M
\mathrm{e}^{\fracd{\EM{{\llvert  z \rrvert }}^2\EM
{{\llvert  a_n \rrvert }}^2}{2}}=\mathrm{e}^{\fracd{\EM{{\llvert  z
\rrvert }}^2s_M^2}{2}}
\end{equation}
with $s_M=\EM{{ (\sum_{n=1}^M \EM{{\llvert  a_n \rrvert }}^2 )}}^{1/2}\le\EM{{\llVert  a\rrVert }}_{{\ell
^2}}$. Hence, for any
$M\in\dN\backslash\EM{{ \{0 \}}}$, $\sum_{n=1}^M a_n
g_n$ is sub-Gaussian
with parameter $s_M$.
Since $a\ne0$, for $M$ large enough, $s_M\ne0$
and then applying Lemma~\ref{DecSG},
\[
\forall t>0,\qquad \dP\EM{{ \Biggl(\sup_{1\le P\le M} \EM{{\Biggl
\llvert \sum_{n=1}^P a_n
g_n \Biggr\rrvert }}> t\EM{{\llVert a\rrVert }}_{{\ell^2}}
\Biggr)}} \le 8\exp\EM{{ \biggl(-\frac{t^2\EM{{\llVert  a\rrVert }}_{{\ell
^2}}^2}{8s_M^2} \biggr)}} \le8
\mathrm{e}^{-\fraca{t^2}{8}}.
\]
Assertion 1 follows
letting $M\to+\infty$.

Let us now prove assertion 2. If there exists $N \in\dN
\backslash\EM{{ \{0 \}}}$, such that
\[
\forall n\ge N,\qquad a_n=0,
\]
then, assertion 2 is fulfilled since
$\sum_{n=1}^{+\infty} a_ng_n=\sum_{n=1}^{N} a_ng_n$ is a
sub-Gaussian random variable with parameter
$s_N=\EM{{ (\sum_{n=1}^N \EM{{\llvert  a_n \rrvert }}^2 )}}^{1/2}
=\EM{{\llVert  a\rrVert }}_{\ell^2}$.
Therefore to prove assertion 2, we now assume that
\[
\forall N\in\dN\backslash\EM{{ \{0 \}}}, \exists n\ge N,\qquad a_n
\ne0,
\]
so that $\sum_{n=N}^{+\infty}\EM{{\llvert  a_n \rrvert }}^2\ne0$ for
any integer $N\ge
1$. Then, applying assertion 1 replacing $a_n$ by $a_n {\mathbf
{1}}_{n\ge N}$, we have
\[
\forall\veps>0, \forall N\in\dN\backslash\EM{{ \{0 \}}},\qquad\dP\EM {{ \Biggl(
\sup_{P\ge N}\EM{{\Biggl\llvert \sum_{n=N}^Pa_ng_n
\Biggr\rrvert }} > \veps \Biggr)}}\le8\mathrm{e}^{-\fraca{\veps^2}{8
{\sum_{n=N}^{+\infty}\EM{{\llvert  a_n \rrvert }}^2}}}.
\]
Since $\EM{{\llVert  a\rrVert }}_{{\ell^2}}^2=\sum_{n=1}^{+\infty
}\EM{{\llvert  a_n \rrvert }}^2<+\infty$, this implies that $\EM
{{ (\sum_{n=1}^N a_ng_n )}}_{N}$
is a Cauchy sequence in probability. Then, by Lemma~3.6 in \cite
{Kallenberg}, the series $\sum_{n=1}^{+\infty} a_ng_n$ converges in
probability. By It\^o--Nisio theorem (see, \cite{LedTal}, e.g.),
this series also converges almost surely, since the random variables
$g_n$, $n\ge1$, are independent.
Moreover, since
$\sup_{M\ge1} s_M^2=\EM{{\llVert  a\rrVert }}_{\ell_2}^2<+\infty$,
equation \eqref{SGSF} implies the uniform integrability of the
sequence $\EM{{ (\mathrm{e}^{\Re\EM{{ (\overline{z} \sum
_{n=1}^M a_n g_n )}}} )}}_{M\ge1}$ for any $z\in\dC$.
Then, letting $M\to+\infty$ in
\eqref{SGSF}, we obtain that $\sum_{n=1}^{+\infty} a_ng_n$ is
sub-Gaussian with parameter $\EM{{\llVert  a\rrVert }}_{{\ell^2}}$.
Moreover,
we conclude the proof noting that
\[
\forall t>0,\qquad \dP\EM{{ \Biggl(\EM{{\Biggl\llvert \sum
_{n=1}^{+\infty} a_ng_n \Biggr
\rrvert }} > t\EM{{\llVert a\rrVert }}_{{\ell ^2}} \Biggr)}} \le\dP\EM {{
\Biggl(\sup_{P\ge1}\EM{{\Biggl\llvert \sum
_{n=1}^Pa_ng_n \Biggr
\rrvert }} > t\EM{{\llVert a\rrVert }}_{{\ell^2}} \Biggr)}}\le8\mathrm
{e}^{-\fraca{t^2}{8}}.
\]
\upqed\end{pf*}
%
%s7 #&#
\section{Main results on conditionally sub-Gaussian series}\label{appB}

%s7.1 #&#
\subsection{Local modulus of continuity}
\label{LMC-A}

\setcounter{equation}{22}

This section is devoted to the proofs of the results stated in
Section~\ref{LMC}.

\begin{pf*}{Proof of Theorem~\ref{CVUS1}}
Let us recall that $x_0\in K_d=\prod_{j=1}^d [a_j,b_j]\subset\dR^d$.
We assume, without loss of generality, that
\[
\forall1\le j\le d,\qquad a_j<b_j.
\]
Actually, if some $a_j=b_j$, we may identify
$\EM{{ (S_N )}}_{N\in\dN}$ and its limits $S$ as random
fields defined on
$K_{d'}\subset\dR^{d'}$ for $d'<d$. Note that if $a_j=b_j$ for all
$1\le j\le d$, there is nothing to prove.

We also assume that $\ga(\omega) \in(0,1)$, which is not restrictive
and allows us to apply equation \eqref{controlrho} as soon as $\EM
{{\llVert  x-y\rrVert }}\le\ga(\omega)$ (with $c_{{2,1}}$ and
$c_{{2,2}}$ which do
not depend on $\ga$).

\textit{First step}. We first introduce a convenient
sequence $\EM{{ (\cD_{\nu_k} )}}_{k\ge1}$ of countable
sets included on dyadics,
which is linked to the quasi-metric $\rho$. It allows to follow some
arguments of the proof of the Kolmogorov's lemma to obtain an upper
bound for the modulus of continuity of $S$.

Let us first introduce some notation. For any $k\in\dN\backslash\EM
{{ \{0 \}}}$ and $j\in\dZ^d$, we set
\[
x_{{k,j}}=\frac{j}{2^k},\qquad \mathcal{D}_k=\EM{{
\bigl\{x_{{k,j}}\dvtx j\in\dZ^d \bigr\}}} \quad\mbox{and}\quad
\nu_k=\min\EM{{ \bigl\{n\in\dN\backslash\EM{{ \{0 \} }}\dvtx
c_{{2,2}} d^{\underline{H}/2}2^{-n\underline{H}}\le 2^{-k} \bigr\}}}
\]
with $c_{{2,2}}$ the constant given by equation \eqref{controlrho}.
Then, choosing $c_{{2,2}}$ large enough (which is not restrictive),
one
checks
that $\EM{{ (\nu_k )}}_{k\ge1}$
is an increasing sequence.
In particular, the sequence $\EM{{ (\cD_{\nu_k} )}}_{k\ge
1}$ is
increasing and
$\cD=\bigcup_{k=1}^{+\infty}\cD_{k}=\bigcup_{k=1}^{+\infty}\cD
_{\nu_k}$. Moreover,
$\cD\cap K_d$ is dense in $K_d$ since $a_j<b_j$ for any $1\le j\le d$.
Then, as done in Step 1 of the proof of Theorem~5.1 of \cite{BL09},
one also checks that for $k$ large enough, $\cD_{\nu_k}\cap K_d$ is a
$2^{-k}$ net of $K_d$ for $\rho$, which means that for any $x\in K_d$,
there exists $j\in\dZ^d$ such that
$
\rho\EM{{ (x,x_{\nu_k,j} )}}\le2^{-k}$,
with $x_{\nu_k,j}=j/2^{\nu_k}\in K_d$.

\textit{Second step}. This step is inspired from Step 2 of \cite
{BL09,BLS11}. The main difference is that we use Proposition~\ref{SGP}
to obtain a uniform control in $N$.

For $k\in\dN\backslash\{0\}$ and $(i,j)\in\mathbb{Z}^d$, we consider
\[
E_{i,j}^{k}= \EM{{ \Bigl\{\omega\dvtx \sup
_{N\in\dN} \bigl\llvert S_{N}\EM{{ (x_{\nu_k,i}
)}}-S_{N}\EM{{ (x_{\nu_k,j} )}}\bigr\rrvert > s\EM{{
(x_{\nu_k,i},x_{\nu_k,j} )}} \varphi\EM {{ \bigl(\rho\EM{{
(x_{\nu _k,i},x_{\nu_k,j} )}} \bigr)}} \Bigr\}}}
\]
with, following \cite{Konog},
%
%e23 #&#
%
\begin{equation}
\label{PHI} \varphi\EM{{ (t )}}=\sqrt{8Ad\log\frac{1}{t}},\qquad t>0,
\end{equation}
for $A>0$ conveniently chosen later.
We choose $\delta\in (0,1 )$ and
set for $k\in\mathbb{N}\backslash\{0\}$,
%
%e24 #&#
%
\begin{eqnarray}
\label{deIk} \delta_k&=&2^{-(1-\delta)k},
\nonumber
\\[-8pt]
\\[-8pt]
I_{{k}}&=&\EM{{ \bigl\{(i,j)\in\EM{{ \bigl(\dZ^d\cap
2^{\nu_k} K_d \bigr)}}^2 \dvtx \rho
(x_{\nu_k,i},x_{\nu_k,j} )\le \delta_k \bigr\}}} \quad
\mbox{and}\quad E_k=\bigcup_{(i,j)\in I_k}
E_{i,j}^{k}.
\nonumber
\end{eqnarray}
Since $\varphi$ is a decreasing function and $s\ge0$, for any $k\in
\dN\backslash\EM{{ \{0 \}}}$ and for any $(i,j)\in I_k$
\[
\dP\EM{{ \bigl(E_{i,j}^{k} \bigr)}}\le\dP\EM{{ \Bigl(\sup
_{N\in
\dN}\bigl\llvert S_{N}\EM{{ (x_{\nu_k,i}
)}}-S_{N}\EM{{ (x_{\nu _k,j} )}}\bigr\rrvert > s\EM{{
(x_{\nu_k,i},x_{\nu
_k,j} )}} \varphi\EM{{ (\delta_k )}}
\Bigr)}}.
\]
Since
$(g_n)_{n\ge1}$ is a sequence of symmetric independent sub-Gaussian
random variables with parameter $s=1$, conditioning to $\EM{{
(W_n )}}_{n\ge1}$ and applying assertion 1 of Proposition~\ref{SGP},
one has
\[
\forall k\in\dN\backslash\EM{{ \{0 \}}}, \forall (i,j)\in I_k,
\qquad \dP\EM{{ \bigl(E_{i,j}^{k} \bigr)}}\le8
\mathrm{e}^{- \fraca
{\varphi\EM{{ (\delta _k )}}^2}{8}} = 8\mathrm {e}^{-Ad(1-\delta) k \log2}
\]
by definition of $s$, $S_N$, $\varphi$ and $\delta_k$.
Moreover, since $K_d\subset\dR^d$ is a compact set, using equation
\eqref{controlrho2} and the definition of $\nu_k$,
one easily proves
that there exists a finite positive constant $c_1\in(0,+\infty)$ such
that for any $k\in\dN\backslash\EM{{ \{0 \}}}$,
$\operatorname{card} I_{{k}}\le
c_{1} 2^{\fracd{2kd}{\underline{H}}}\delta_k^{\fraca{d}{\overline{H}}}$.
Hence,
\[
\sum_{k=1}^{+\infty}\dP\EM{{ (E_k
)}}\le \sum_{k=1}^{+\infty}\sum
_{\EM{{ (i,j )}}\in
I_{{k}}}\mathbb{P} \bigl(E_{i,j}^{k} \bigr)
\le{c_{1}} \sum_{k=1}^{+\infty}
\mathrm{e}^{-\EM{{ (A(1-\delta) -\fraca
{2}{\underline{H}}+\fracd{1-\delta}{\overline{H}} )}} k d\log2} <+\infty
\]
choosing $A>\frac{2}{\underline{H}}-\frac{1}{\overline{H}}$ and
$\delta$ small enough.
Then, setting
\[
\Omega''=\Omega'\cap\EM{{
\Biggl(\bigcup_{k=1}^{+\infty} \bigcap
_{\ell=k}^{+\infty} E_\ell^c \Biggr)}}
\]
with $\Omega'$ the almost sure event introduced by Assumption~\ref{Hx0},
the Borel--Cantelli lemma leads to
$\mathbb{P}(\Omega'')=1$. Moreover, by Assumption~\ref{Hx0},
for any $\omega\in\Omega''$
there exists $k^*\EM{{ (\omega )}}$ such that for every
$k\ge k^*\EM{{ (\omega )}}$ and for all $x, y \in
{\mathcal D}_{\nu_k}$ with $x,y \in
B(x_0,\gamma(\omega))\cap K_d$
and $\rho\EM{{ (x,y )}}\le\delta_k=2^{-(1-\delta)k}$,
%
%e25 #&#
%
\begin{equation}
\label{maj:alpha:rho} \sup_{N\in\dN}\EM{{\bigl\llvert S_{N}
\EM{{ (x )}}-S_{N}\EM {{ (y )}} \bigr\rrvert }}\le C \rho
(x,y)^{\beta} \bigl |\log\bigl( \rho(x,y)\bigr)\bigr |^{\eta+1/2}.
\end{equation}

\textit{Third step}.
In this step, we prove that \eqref{maj:alpha:rho} holds, up to a
multiplicative constant, for any $x,y\in\cD$ closed enough to $x_0$.
This step is adapted from Step 4 of the proof of Theorem~5.1 in \cite
{BL09}, taking care that \eqref{maj:alpha:rho} only holds for some
$x,y\in\cD_{\nu_k}\cap K_d$ randomly closed enough of $x_0$. Let us
mention that this step has been omitted in the proof of the main result
of \cite{BLS11} but is not trivial. We then decide to provide a proof
here for the sake of completeness and clearness.

Let us now fix $\omega\in\Omega''$ and denote by $\kappa\ge1$ the
constant appearing in the quasi-triangle inequality satisfied by $\rho
$. We also consider the function $F$ defined on $(0,+\infty)$ by\vspace*{-1pt}
\[
F(h):=h^{\beta} \bigl |\log( h)\bigr |^{\eta+1/2}.
\]
Observe that $F$ is a random function since $\beta$ and $\eta$ are
random variables.
Then, we choose $k_0=k_0(\omega)\in\dN$ such that the three
following assertions are fulfilled:
\renewcommand\theenumi{(\alph{enumi})}\renewcommand\labelenumi{\theenumi}{
\begin{enumerate}[(a)]
\item
$F$ is increasing on $(0,\delta_{k_0}]$, where $\delta_{k}$ is given
by \eqref{deIk},
\item for all $k\ge k_0(\omega)$, ${\mathcal D}_{\nu_k}\cap K_d$ is
a $2^{-k}$ net of $K_d$ for $\rho$,
\item
$ 2^{k_0}\delta_{k_0+1}>3\kappa^2$.
\end{enumerate}
}
Even if it means to choose $k^*\EM{{ (\omega )}}$ larger,
we can assume
that $k^*\EM{{ (\omega )}}\ge k_0$ and that
%
%e26 #&#
%
\begin{equation}
\label{Cga} \gamma(\omega)\ge \biggl(\frac{\delta_{k^*(\omega)}}{3\kappa^2 c_{{2,2}}} \biggr)^{1/\underline{H}}
:=2\gamma^*(\omega),
\end{equation}
where $\underline{H}$ and $c_{{2,2}}$ are defined in equation \eqref
{controlrho}.\vadjust{\goodbreak}

Let us now consider $x,y \in\mathcal{D}\cap K_d\cap B(x_0,\ga
^*(\omega)) $ such that $x\neq y$. Let us first note that $x,y\in
B(x_0,\ga(\omega))$. Moreover, since $\EM{{\llVert  x-y\rrVert
}}\le2\ga^*(\omega
)\le\gamma(\omega)\le1$,
the upper bound of equation \eqref{controlrho} leads to
\[
3\kappa^2\rho(x,y)\le3\kappa^2 c_{{2,2}}\EM{{
\llVert x-y\rrVert }}^{\underline
{H}}\le\delta_{k^*(\omega)}
\]
by definition of $\ga^*(\omega)$.
Then,
there exists a unique $k\ge k^*(\omega)$ such that
%
%e27 #&#
%
\begin{equation}
\label{choixk} \delta_{k+1}<3\kappa^2\rho(x,y)\le
\delta_k.
\end{equation}
Furthermore, since $x, y \in\mathcal{D}\cap K_d$, there exists $n\ge
k+1$ such that $x,y\in\cD_{\nu_n}\cap K_d$ and for $j=k,\ldots
,n-1$, there exist $x^{(j)}\in\cD_{\nu_j}\cap K_d$ and $ y^{(j)}\in
\cD_{\nu_j}\cap K_d$ such that
%
%e28 #&#
%
\begin{equation}
\label{refIn} \rho\EM{{ \bigl(x,x^{(j)} \bigr)}}\le2^{-j}
\quad\mbox{and}\quad \rho\EM{{ \bigl(y,y^{(j)} \bigr)}}
\le2^{-j}.
\end{equation}
Let us now fix
$N\in\dN$ and focus on $S_N(x)-S_N(y)$. Then, setting $x^{(n)}=x$ and
$y^{(n)}=y$,
%
%e29 #&#
%
\begin{eqnarray}
\label{Zdec} S_N(x)-S_N(y)&=& \EM{{
\bigl(S_N\EM{{ \bigl(x^{(k)} \bigr)}}-S_{N}\EM{{
\bigl(y^{(k)} \bigr)}} \bigr)}}+ \sum_{j=k}^{n-1}
\EM{{ \bigl(S_{N}\EM{{ \bigl(x^{(j+1)} \bigr)}}-S_{N}
\EM{{ \bigl(x^{(j)} \bigr)}} \bigr)}}
\nonumber
\\[-8pt]
\\[-8pt]
&&{} -\sum_{j=k}^{n-1}\EM{{
\bigl(S_{N}\EM{{ \bigl(y^{(j+1)} \bigr)}}-S_{N}\EM{{
\bigl(y^{(j)} \bigr)}} \bigr)}}.
\nonumber
\end{eqnarray}

The following lemma, whose proof is {given below}
for the sake of clearness, allows
to apply \eqref{maj:alpha:rho} for each term of the right-hand side of
the last inequality.

%le7.1 #&#
%
\begin{lemma}
\label{LemCS} Choosing $k^*(\omega)$ large enough,
the sequences $\EM{{ (x^{(j)} )}}_{k\le j\le n}$ and $\EM
{{ (y^{(j)} )}}_{k\le
j\le n}$ satisfy the three following assertions.
\begin{enumerate}[3.]
\item$x^{(j)},y^{(j)}\in B(x_0,\gamma(\omega))$ for any $j=k,\ldots
, n$,
\item for any $j=k,\ldots, n-1$, $\max(\rho(x^{(j+1)},x^{(j)}),\rho
(y^{(j+1)},y^{(j)}))\le\delta_{j+1}$,
\item$\rho(x^{(k)},y^{(k)})\le\delta_k$.
\end{enumerate}
\end{lemma}

Therefore, even if it means to choose $k^*\EM{{ (\omega
)}}$ larger,
applying this lemma and equations \eqref{maj:alpha:rho} and \eqref
{Zdec}, we obtain
\[
\EM{{\bigl\llvert S_{N}\EM{{ (x )}}-S_{N}\EM{{ (y )}}
\bigr\rrvert }}\le C\EM{{ \Biggl(F \bigl(\rho\bigl(x^{(k)},y^{(k)}
\bigr) \bigr) + 2\sum_{j=k}^{n-1} F\EM{{ (
\delta_{j+1} )}} \Biggr)}}
\]
since $F$ is increasing on $(0,\delta_{k_0}]$ and since $j\ge k_0$.
This implies, by definition of $F$ that
\[
\EM{{\bigl\llvert S_{N}\EM{{ (x )}}-S_{N}\EM{{ (y )}}
\bigr\rrvert }}\le C\EM{{ \bigl(F \bigl(\rho\bigl(x^{(k)},y^{(k)}
\bigr) \bigr) + 2\tilde{C} F\EM{{ (\delta_{k+1} )}} \bigr)}},
\]
where
\[
\tilde{C}(\omega)=2\sum_{j=0}^{+\infty}
\delta_j^{\beta(\omega)} (j+1)^{\max(\eta(\omega)+1/2,0)}<+\infty
\]
since $\beta>0$ and $\delta_j=2^{-(1-\delta)j}$ with $\delta<1$.
Then, since $F$ is increasing on $(0,\delta_0)$, by assertion 3
of Lemma~\ref{LemCS} and equation \eqref{choixk}, we get
\[
\EM{{\bigl\llvert S_{N}\EM{{ (x )}}-S_{N}\EM{{ (y )}}
\bigr\rrvert }}\le C(1+2\tilde{C})F\bigl(3\kappa^2\rho(x,y)\bigr),
\]
for every
$N\in\dN$ and $x,y\in\cD\cap B(x_0, \ga^*(\omega))\cap K_d$.
Therefore, by continuity of $\rho$ and each $S_N$ and by density of
$\cD\cap K_d$ in $K_d$
%
%e30 #&#
%
\begin{equation}
\label{ModSNC} \EM{{\bigl\llvert S_{N}\EM{{ (x )}}-S_{N}
\EM{{ (y )}} \bigr\rrvert }}\le C(1+2\tilde{C})F\bigl(3\kappa^2
\rho(x,y)\bigr),
\end{equation}
for every
$N\in\dN$ and $x,y\in B(x_0, \ga^*(\omega))\cap K_d$.

\textit{Fourth step: Uniform convergence of ${S}_{N}$}.
Let us now consider
\[
\tilde{\Omega}=\bigcap_{u\in\cD}\EM{{ \Bigl\{\lim
_{N\to+\infty
} S_N(u)= S(u) \Bigr\}}}\cap
\Omega''.
\]
Observe that $\dP\EM{{ (\tilde{\Omega} )}}=1$. Let us
now fix $\omega
\in\tilde{\Omega}$. Hence, by equation \eqref{ModSNC}, the sequence
$ \EM{{ (S_{N}\EM{{ (\cdot )}}(\omega)
)}}_{N\in\dN}$, which converges
pointwise on $\cD\cap B(x_0, \ga^*(\omega))\cap K_d$ is uniformly
equicontinuous on $B(x_0, \ga^*(\omega))$. Since $\cD\cap B(x_0, \ga
^*(\omega))\cap K_d$ is dense in $B(x_0, \ga^*(\omega))\cap K_d$, by
Theorem I.26 and adapting Theorem I.27 in \cite{Reed}, $\EM{{
(S_{N}\EM{{ (\cdot )}}(\omega) )}}_{N\in\dN}$
converges uniformly on
$B(x_0, \ga^*(\omega))\cap K_d$. Therefore, its limit $S$ is
continuous on $B(x_0, \ga^*(\omega))\cap K_d$. Moreover, letting
$N\to+\infty$ in \eqref{ModSNC} (which holds since $\omega\in
\tilde{\Omega}$), we get
%
%e31 #&#
%
\begin{equation}
\label{ModSC} \EM{{\bigl\llvert S\EM{{ (x )}}-S\EM{{ (y )}} \bigr\rrvert }}\le
C(1+2\tilde{C})F\bigl(3\kappa^2\rho(x,y)\bigr),
\end{equation}
for every $x,y\in B(x_0, \ga^*(\omega))\cap K_d$, which concludes the proof.
\end{pf*}

Let us now prove Lemma~\ref{LemCS}.
\begin{pf*}{Proof of Lemma~\ref{LemCS}}
Let us first observe that $x^{(n)}=x\in B(x_0,\gamma(\omega))\cap
K_d$ and $y^{(n)}=y\in B(x_0,\gamma(\omega))\cap K_d$. Let us now fix
$j\in\EM{{ \{k,\ldots,n-1 \}}}$.
The lower bound of equation \eqref{controlrho} leads to
\[
\EM{{\bigl\llVert x^{(j)}-x_0\bigr\rrVert }}\le\EM{{\bigl
\llVert x^{(j)}-x\bigr\rrVert }}+\EM{{\llVert x-x_0\rrVert }}
\le\frac{\rho\EM{{
(x^{(j)},x )}}^{1/\overline{H}}}{c_{{2,1}}^{1/\overline
{H}}}+\EM{{\llVert x-x_0\rrVert }}.
\]
Since $x\in B(x_0,\gamma^*(\omega))$ with $\gamma^*$
satisfying equation \eqref{Cga} and since $\rho\EM{{
(x^{(j)},x )}}\le
2^{-j}$ with $j\ge k\ge k^*(\omega)$,
\[
\EM{{\bigl\llVert x^{(j)}-x_0\bigr\rrVert }}\le
\frac{2^{-k^*(\omega
)/\overline
{H}}}{c_{{2,1}}^{1/\overline{H}}}+\frac{\gamma(\omega)}{2}.
\]
Then, choosing $k^*(\omega)$ large enough,
$x^{(j)}\in B(x_0,\gamma(\omega))$ for $j=k,\ldots,n-1$. The same
holds for $y^{(j)}$. Assertion 1 is then proved.

Let us now observe that since $j\ge k_0$ and since $\kappa\ge1$,
%
%e32 #&#
%
\begin{equation}
\label{eqB} 2^j\delta_{j+1}\ge2^{k_0}
\delta_{k_0+1}>3\kappa^2\ge3\kappa
\end{equation}
by definition of $k_0$ (see the third step of the proof of Theorem~\ref
{CVUS1}).
Then, using the quasi-triangle inequality fulfilled by $\rho$ and
\eqref{refIn}, we obtain that
\[
\rho\bigl(x^{(j+1)},x^{(j)}\bigr)\le3\kappa2^{-(j+1)}\le
\frac{\delta
_{j+1}}{2}\le\delta_{j+1}.
\]
Since the same holds for $\rho(y^{(j+1)},y^{(j)})$, assertion 2 is fulfilled.
Moreover, applying twice the quasi-triangle inequality fulfilled by
$\rho$ and equations \eqref{choixk}, \eqref{refIn} and \eqref{eqB}
(with $j=k$), we obtain
\[
\rho\bigl(x^{(k)},y^{(k)}\bigr)\le\kappa^2
\bigl(2^{1-k}+\rho(x,y)\bigr)\le3\kappa ^2\rho(x,y) \le
\delta_k,
\]
which is assertion 3.
\end{pf*}

Let us now focus on Corollary~\ref{ArgReco}. Its proof is based on the
following technical lemma.

%le7.2 #&#
%
\begin{lemma}
\label{LemReco}
Let $K_d=\prod_{i=1}^d[a_i,b_i]\subset\dR^d$ be a compact
$d$-dimensional interval, $\beta\in(0,1)$, $\eta\in\dR$ and $\rho
$ be a quasi-metric on $\dR^d$ satisfying equation \eqref{controlrho}.
Let $\EM{{ (f_n )}}_{n\in\dN}$ be sequence of functions
defined on $K_d$ and
let $\EM{{ (\mathring{B}(x_i,r_i) )}}_{1\le i\le p}$ be a
finite covering
of $K_d$ by open balls with $x_i\in K_d$ and $r_i>0$. Assume that for
each $1\le i\le p$, there exists a finite positive constant $C_i$ such that
\[
\forall x,y\in\mathring{B}(x_i,r_i)\cap
K_d,\qquad \sup_{n\in\dN}\EM {{\bigl\llvert
f_n(x)-f_n(y) \bigr\rrvert }} \le C_i
\rho\EM{{ (x,y )}}^\beta\EM{{ \bigl[\log\EM{{ \bigl(1+\rho
(x,y)^{-1} \bigr)}} \bigr]}}^\eta.
\]
Then there exists a finite positive constant $C$ such that
%
%e33 #&#
%
\begin{equation}
\label{CUf} \forall x,y\in K_d,\qquad\sup_{n\in\dN}
\EM{{\bigl\llvert f_n(x)-f_n(y) \bigr\rrvert }} \le C
\rho \EM{{ (x,y )}}^\beta\EM{{ \bigl[\log\EM{{ \bigl(1+\rho
(x,y)^{-1} \bigr)}} \bigr]}}^\eta.
\end{equation}
\end{lemma}
\begin{pf} By the Lebesgue's number lemma, there exists $r>0$ such that
\[
\forall x\in K_d, \exists1\le i\le p,\qquad\mathring{B}(x,r) \subset
\mathring{B}\EM{{ (x_i,r_i )}}.
\]
Let us first note that since the map $F_{\rho} \dvtx(u,v)\mapsto\rho
^{\beta}(u,v)\EM{{ [\log\EM{{ (1+\rho(u,v)^{-1}
)}} ]}}^\eta$ is positive
and continuous on the compact set $\tilde{K}=\EM{{ \{(u,v)\in
K\times K / \EM{{\llVert  u-v\rrVert }}\ge r \}}}$,
\[
m:=\inf_{\tilde{K}} F_{\rho} \in(0,+\infty).
\]
Then distinguishing the cases $\EM{{\llVert  x-y\rrVert }}< r$ and
$\EM{{\llVert  x-y\rrVert }}\ge r$,
one easily sees that
\[
\sup_{n\in\dN} \EM{{\bigl\llvert f_n(x)-f_n(y)
\bigr\rrvert }}\le\max\EM {{ \biggl(\max_{1\le i\le p}
C_i, \frac{M}{m} \biggr)}} \rho (x,y)^{\beta}\EM{{
\bigl[\log\bigl(1+\rho (x,y)^{-1}\bigr) \bigr]}}^{\eta},
\]
where $M=\sup_{n\in\dN}\sup_{x,y\in K_d} \EM{{\llvert
f_n(x)-f_n(y) \rrvert }}$.
It remains to prove
that $M<+\infty$.
Note that
\[
\mathop{\sup_{x,y\in K_d}}_{\EM{{\llVert  x-y\rrVert }}<r}\EM {{\bigl\llvert
f_n(x)-f_n(y) \bigr\rrvert }}\le c \max
_{1\le i\le p} C_i,
\]
where $c=\sup_{K_d\times K_d}F_{\rho}
<+\infty$ by continuity of $ F_{\rho}$
on the compact set $K_d\times K_d$. Then since $K_d$ is a compact
convex set, using a chaining argument, one easily obtains that
$M<+\infty$, which concludes the proof.
\end{pf}

\begin{pf*}{Proof of Corollary~\ref{ArgReco}}
We only prove
assertion 1. Actually, assertion 2 is proved using
the same arguments but replacing $K_d$ by $B(x_0,\gamma)\cap K_d$.

Assume that for any $x_0\in K_d$, Assumption~\ref{Hx0} holds with
$\Omega'$, $\beta$, $\eta$ and the quasi-metric $\rho$ independent
of $x_0$. Following the proof of
Theorem~\ref{CVUS1} and keeping its notation, let us quote that
$\gamma^*$ and $\tilde{\Omega}$
do not depend on $x_0$.
Let us now fix $\omega\in\tilde{\Omega}$.
From the third step of the proof of Theorem~\ref{CVUS1} and Lemma~\ref
{LemReco}, we deduce that equation \eqref{ModSNC} still holds for any
$x,y\in K_d$. This allows to replace $B\EM{{ (x_0,\gamma
^*(\omega) )}}$
by $K_d$ in the fourth step of the proof of Theorem~\ref{CVUS1}, which
leads to assertion 1.
\end{pf*}

%s7.2 #&#
\subsection{Rate of almost sure uniform convergence}
\label{CVUS2-A}

\begin{pf*}{Proof of Theorem~\ref{CVUS2}}
Let us first observe that Theorem~\ref{CVUS1} holds. Then, for almost
$\omega$, even if it means to choose $\ga$ smaller, the sequence of
continuous functions $\EM{{ (S_N(\cdot)(\omega) )}}_{N\in
\dN}$
converges uniformly on $B(x_0,\ga(\omega))\cap K_d$, which implies
that each $R_N(\cdot)(\omega)$ is continuous on $B(x_0,\ga(\omega
))\cap K_d$. As in the proof of Theorem~\ref{CVUS1}, we assume without
loss of generality that $K_d=\prod_{j=1}^d[a_j,b_j]$ with $a_j<b_j$.

\begin{pf*}{Proof of assertion 1} Since it is
quite similar to the proof of equation \eqref{ModSNC},
we only sketch it.

For $k\in\dN\backslash\{0\}$, $N\in\dN$ and $(i,j)\in\mathbb
{Z}^d$, we consider
\[
E_{i,j}^{k,N}= \EM{{ \bigl\{\omega\dvtx\bigl\llvert
R_{N}\EM{{ (x_{\nu_k,i} )}}-R_{N}\EM{{
(x_{\nu _k,j} )}}\bigr\rrvert > \sqrt{\log (N+2)} r_N\EM{{
(x_{\nu_k,i},x_{\nu_k,j} )}} \varphi\EM {{ \bigl(\rho\EM{{
(x_{\nu _k,i},x_{\nu_k,j} )}} \bigr)}} \bigr\}}}
\]
with
$r_N$ defined by \eqref{rN}, $\varphi$ by \eqref{PHI} and $\EM
{{ (\nu _k )}}_{k\ge1}$ by Step 1 of the proof of Theorem~\ref{CVUS1}. Then, we
proceed as in Step 1 of
the proof of Theorem~\ref{CVUS1} replacing the set $E_k$ by
\[
E_k'=\bigcup_{N=0}^{+\infty}
\bigcup_{(i,j)\in I_k} E_{i,j}^{k,N},
\]
with $I_k$ and $\delta_k$ defined by \eqref{deIk},
and applying assertion 2 of Proposition~\ref{SGP} instead of
assertion~1. Then,
choosing the constant $A$, which appears in the definition of $\varphi
$, and $\delta\in(0,1)$ such that
\[
A(1-\delta) -\frac{2}{\underline{H}}+\frac{1-\delta}{\overline
{H}} >0 \quad\mbox{and}\quad A(1-
\delta)\log2>1
\]
we obtain that
\[
\sum_{k=1}^{+\infty}\dP\EM{{
\bigl(E_k' \bigr)}}\le c_2 \sum
_{N=2}^{+\infty} 2^{-A(1-\delta)\log N }=c_2\sum
_{N=2}^{+\infty
} N^{-A(1-\delta)\log2 }<+\infty
\]
with $c_{2}$ a finite positive constant.
Then, by Borel--Cantelli lemma, the definition of $\varphi$ and
Assumption~\ref{Hx02}, almost surely there exists an integer $ k^*\EM
{{ (\omega )}}$ such that for every $k\ge k^*\EM{{
(\omega )}}$, for all $N\in
\dN$, and for all $x, y \in{\mathcal D}_{\nu_k}$ with $x,y \in
B(x_0,\gamma(\omega))\cap K_d$ and $\rho\EM{{ (x,y
)}}\le\delta
_k=2^{-(1-\delta)k}$
\[
\EM{{\bigl\llvert R_{N}\EM{{ (x )}}-R_{N}\EM{{ (y )}}
\bigr\rrvert }}\le C b(N)\sqrt{\log(N+2)}\rho (x,y)^{\beta} \bigl |\log\bigl(
\rho(x,y)\bigr)\bigr |^{\eta+1/2}.
\]
In addition, replacing in Step 2 of the proof of Theorem~\ref{CVUS1},
$S_N$ by $R_N$ (which still be, for almost all $\omega$, continuous on
$B(x_0,\ga(\omega))\cap K_d$), we obtain that for almost all $\omega$,
there exists $\ga^*\in(0,\ga)$, such that
%
%e34 #&#
%
\begin{equation}
\label{modR} \EM{{\bigl\llvert R_{N}\EM{{ (x )}}-R_{N}
\EM{{ (y )}} \bigr\rrvert }}\le C b(N)\sqrt{\log(N+2)}\rho (x,y)^{\beta} \bigl |
\log\bigl( \rho(x,y)\bigr)\bigr |^{\eta+1/2}
\end{equation}
for every $N\in\dN$ and $x,y\in B(x_0, \ga^*(\omega))\cap K_d$.
This establishes assertion 1.
\noqed
\end{pf*}

\begin{pf*}{Proof of assertion 2}
This assertion
follows from equations \eqref{modR} and \eqref{Contx0},
the continuity of $\rho$ on the compact set $B(x_0,\ga(\omega))\cap
K_d$ and
\[
\EM{{\bigl\llvert R_N(x) \bigr\rrvert }}\le\EM{{\bigl\llvert
R_N(x)-R_N(x_0) \bigr\rrvert }}+\EM{{\bigl
\llvert R_N(x_0) \bigr\rrvert }}.
\]
The proof of Theorem~\ref{CVUS2} is then complete.\qed
\noqed
\end{pf*}
\noqed
\end{pf*}

\begin{pf*}{Proof of Corollary~\ref{ArgReco2}}
We only prove
assertion 1. Actually, assertion 2 is proved using
the same arguments but replacing $K_d$ by $B(x_0,\gamma)\cap K_d$.

Let us assume that Assumption~\ref{Hx02} holds with $\Omega'$, $\beta
$, $\eta$ and the quasi-metric $\rho$ independent of $x_0$.
Note first that the almost sure event $\tilde{\Omega}$ under which
\eqref{modR} holds does not depend on $x_0$. Then applying Lemma~\ref
{LemReco} to $f_n=R_n/\EM{{ (b(n)\sqrt{\log(n+2)} )}}$,
we obtain that
equation \eqref{modR} still holds for $x,y\in K_d$.
If moreover for some $x_0$, equation \eqref{Contx0} is fulfilled, then
following the proof of assertion 2 of Theorem~\ref{CVUS2},
we also have:
there exists $C$ a finite positive random variable such that for all
$N\in\dN$,
\[
\sup_{x\in K_d} \EM{{\bigl\llvert R_N(x) \bigr\rrvert
}}\le C b(N)\sqrt{\log \EM{{ (N+2 )}}},
\]
which concludes the proof.
\end{pf*}
%

%s8 #&#
\section{Shot noise series}
%s8.1 #&#
\subsection{Proof of Theorem \texorpdfstring{\protect\ref{uCOLALE}}{4.1}}
\label{uCOLALE-A}

\setcounter{equation}{34}

Let $\EM{{ (g_n )}}_{n\ge1}$ be a Rademacher sequence,
that is, a sequence
of i.i.d. random variables with symmetric Bernoulli distribution.
This Rademacher sequence is assumed to be independent of $\EM{{
(T_n,X_n )}}_{n\ge1}$. Then, by independence and also by
symmetry of the
sequence $\EM{{ (X_n )}}_{n\ge1}$,
$(X_ng_n)_{n\ge1}$ has the same distribution as $\EM{{
(X_n )}}_{n\ge1}$
and is independent of the sequence $\EM{{ (T_n )}}_{n\ge1}$.

Let us now set
\[
W_n(\al)=T_n^{-1/\al}X_n \quad
\mbox{and}\quad S_N(\al)=\sum_{n=1}^NW_n(
\al)g_n,
\]
so that $\{S_N^*(\al), N\ge1\}$ has the same finite distribution as
$\{S_N(\al), N\ge1\}$. Moreover, since
\[
\sum_{n=1}^N\bigl |W_n(
\al)\bigr |^2=\sum_{n=1}^NT_n^{-2/\al}|X_n|^2
\]
with $X_n\in L^{2p}$ (with $p>0$),
Assumption~\ref{HypTVg} is fulfilled on any $K_1=[a,b]\subset
(0,\min(2,2p))$ (see, e.g., \cite{sam:taqqu}).

Let us now fix $a,b\in(0,\min(2,2p))$ such that $a<b$, $a'\in(0,a)$
and $b'\in(b,\min(2,2p))$.

\begin{pf*}{Proof of assertion 1}
By the Mean Value
Inequality, we get that
for any $\al,\al'\in[a,b]$ and $n\ge1$,
%
%e35 #&#
%
\begin{equation}
\label{majT} \EM{{\bigl\llvert T_n^{-1/\al}-T_n^{-1/\al'}
\bigr\rrvert }}\le c\EM{{\bigl\llvert \al-\al' \bigr\rrvert }}\max
\EM{{ \bigl(T_n^{-1/b'},T_n^{-1/a'}
\bigr)}}
\end{equation}
with $c$ a finite positive constant.
It follows that, almost surely, for all $\alpha,\alpha'\in[a,b]$,
%
%e36 #&#
%
\begin{equation}
\label{senalpha} s\bigl(\alpha,\alpha'\bigr):= \Biggl(\sum
_{n=1}^{+\infty}\bigl |W_n(\al)-W_n
\bigl(\al '\bigr)\bigr |^2 \Biggr)^{1/2}\le C\bigl |\alpha-
\alpha'\bigr |,
\end{equation}
with $C=c (\sum_{n=1}^{+\infty}T_n^{-2/b'}|X_n|^2+ \sum_{n=1}^{+\infty}T_n^{-2/a'}|X_n|^2 )^{1/2}<+\infty$ since
$|X_n|^2\in L^p$ with $2p> b'\ge a'$ and $a',b'\in(0,2)$.
Therefore, the assumptions of assertion 1 of Corollary~\ref{ArgReco} hold.
Let us now remark that for all $\al,\al'\in[a,b]$,
\[
\bigl\{\bigl(S_N^*(\al)-S_N^*\bigl(\al'
\bigr), s\bigl(\al,\al'\bigr)\bigr); N\ge1\bigr\}\stackrel {\mathrm{fdd}} {=}
\bigl\{\bigl(S_N(\al)-S_N\bigl(\al'\bigr), s
\bigl(\al,\al'\bigr)\bigr); N\ge1\bigr\}.
\]
This allows us to replace $S_N$ by $S_N^*$ in the second step of the
proof of Theorem~\ref{CVUS1}. Then, the third and the fourth step of
this proof still hold replacing $S_N$ by $S_N^*$ and the limit $S$ by
the limit $S^*$ since each $S_N^*$ is continuous (as $S_N$ is) and
since $S_N^*$ converges pointwise to $S^*$. This allows us to also
replace $(S_N,S)$ by $(S_N^*,S^*)$ in the proof of assertion 1 of
Corollary~\ref{ArgReco}.
It follows that
almost surely, $\EM{{ (S_N^* )}}_{N\in\dN}$ converges
uniformly on
$[a,b]$ to $S^*$.
Since this holds for any $0<a<b<\min(2,2p)$, assertion 1 of Theorem~\ref{uCOLALE} is established.
\noqed
\end{pf*}

\begin{pf*}{Proof of assertion 2}
Since almost surely the sequence of continuous random fields $\EM
{{ (S_N^* )}}_{N\in\dN}$ converges uniformly on $[a,b]$,
for all $N\in\dN
$ the rest $R_N^*$, defined by
\[
R_N^*(\alpha):=\sum_{n=N+1}^{+\infty}T_n^{-1/\al}X_n,
\]
is also continuous on $[a,b]$.
Remark also that we have, for all $\al,\al'\in[a,b]$ and $N\in\dN$,
\[
\bigl(R_N^*(\al)-R_N^*\bigl(\al'\bigr),
r_N\bigl(\al,\al'\bigr)\bigr)\stackrel{d} {=}
\bigl(R_N(\al )-R_N\bigl(\al'\bigr),
r_N\bigl(\al,\al'\bigr)\bigr),
\]
where $R_N(\al)=\sum_{n=N+1}^{+\infty}T_n^{-1/\al}X_ng_n=S(\al
)-S_N(\al)$ and
\[
r_N\EM{{ \bigl(\al,\al' \bigr)}}=\EM{{ \Biggl(\sum
_{n=N+1}^{+\infty}|X_n|^2
\EM{{\bigl\llvert T_n^{-1/\alpha}-T_n^{-1/\alpha
'}
\bigr\rrvert }}^2 \Biggr)}}^{1/2}.
\]
As done for $S_N$, the previous lines allow to replace $R_N$ by $R_N^*$
in the proof of Theorem~\ref{CVUS2}. Moreover,
by
equation \eqref{majT},
almost surely, for all $N\in\dN$, and $\al,\al'\in[a,b]$,
%
%e37 #&#
%
\begin{equation}
\label{rNenalpha} r_N\EM{{ \bigl(\al,\al' \bigr)}}\le c\bigl |
\al-\al'\bigr | \EM{{ \Biggl(\sum_{n=N+1}^{+\infty}
T_n^{-2/b'} \EM{{\llvert X_n \rrvert
}}^2+\sum_{n=N+1}^{+\infty}T_n^{-2/a'}|X_n|^2
\Biggr)}}^{1/2}.
\end{equation}
Let us now fix $p'>0$ such that $1/p'\in(0,1/b-1/\min(2p,2))$.
Choosing if necessary $b'>b$ smaller, we assume without loss of
generality that $1/p'\in(0,1/b'-1/\min(2p,2))\subset(0,1/a'-1/\min
(2p,2))$. Then, by Theorem~2.2 in \cite{COLALE07}, almost surely, for
all $\al,\al'\in[a,b]$,
\[
\sup_{N\in\dN}{N}^{2/p'} \EM{{ \Biggl(\sum
_{n=N+1}^{+\infty} T_n^{-2/b'} \EM{{
\llvert X_n \rrvert }}^2+\sum
_{n=N+1}^{+\infty
}T_n^{-2/a'}|X_n|^2
\Biggr)}}<+\infty
\]
since $X_n^2\in L^p$ with $p>b'/2>a'/2$ and $a',b'\in(0,2)$.
Note also that by Theorem~2.1 in \cite{COLALE07}, for all $x_0=\al
_0\in[a,b]$, almost surely
\[
\sup_{N\in\dN} {N}^{1/p'}\EM{{\Biggl\llvert \sum
_{n=N+1}^{+\infty} T_n^{-1/\al _0}
X_n \Biggr\rrvert }} <+\infty.
\]
Therefore, the assumptions of assertion 1 of Corollary~\ref{ArgReco} hold
with $b(N)=(N+1)^{-1/p'}$ for any $p'$ such that $1/p'\in(0,1/b-1/\min
(2p,2))$. And then, substituting in its proof $R_N$ by $R_N^*$, almost surely
\[
\sup_{N\in\dN} \sup_{\al\in[a,b]} {N}^{1/p'}
\EM{{\bigl\llvert R_N^*\EM{{ (\al )}} \bigr\rrvert }} <+\infty,
\]
which concludes the proof.
\noqed\end{pf*}

%s8.2 #&#
\subsection{Modulus of continuity and rate of convergence}
\label{SNMC-A}

This section is devoted to the proofs of the results stated in
Section~\ref{MCSN}. First, let us establish Theorem~\ref{CrMom}.
\begin{pf*}{Proof of Theorem~\ref{CrMom}}
Let {us fix} $x_0=\EM{{ (\al_0,u_0 )}}\in
K_{d+1}=[a,b]\times\prod_{j=1}^d[a_j,b_j] \subset(0,2)\times\dR^d
$.

\begin{pf*}{Proof of assertion 1}
Let us assume that $p> b/2$ and consider $s$ the conditional parameter
defined by \eqref{VCS}. Then,
for any $x=(\al,u)\in K_{d+1}$ and $y=(\al',v)\in K_{d+1}$,
%
%e38 #&#
%
\begin{equation}
\label{decs} s\EM{{ (x,y )}}\le s_{1}\EM{{ (x,y )}}+
s_{2}\EM{{ (x,y )}},
\end{equation}
where
\[
s_{1}(x,y)=\EM{{ \Biggl(\sum_{n=1}^{+\infty}
T_n^{-2/\al}\EM{{\bigl\llvert V_n(x)-V_n(y)
\bigr\rrvert }}^2 \Biggr)}}^{1/2}
\]
and
\[
s_{2}(x,y)=\EM{{ \Biggl(\sum_{n=1}^{+\infty}
\EM{{ \bigl(T_n^{-1/\al
}-T_n^{-1/\al '}
\bigr)}}^2 \EM{{\bigl\llvert V_n(y) \bigr\rrvert
}}^2 \Biggr)}}^{1/2}.
\]
First, let us focus on $s_1$. Note that
for any $x,y\in K_{d+1}$,
\[
s_1\EM{{ (x,y )}}\le C_1 \rho\EM{{ (x,y
)}}^{\beta} \log\bigl(1+\rho (x,y)^{-1}\bigr)^\eta,
\]
with $C_1=\EM{{ (\sum_{n=1}^{+\infty} T_n^{-2/b}\EM{{\llvert
Y_n \rrvert }}^2+\sum_{n=1}^{+\infty} T_n^{-2/a}\EM{{\llvert  Y_n
\rrvert }}^2 )}}^{1/2}$, where we have set
\[
Y_n= \mathop{\sup_{x,y\in K_{d+1}}}_{x\ne y}
\frac{\EM{{\llvert
V_n(x)-V_n(y) \rrvert }}}{\rho(x,y)^{\beta}\log(1+\rho
(x,y)^{-1})^\eta}.
\]
Since $K_d$ is a convex compact set, applying a chaining argument and
using the continuity of $\rho$, one checks that equation \eqref{UBV}
implies that $Y_n\in L^{2p}$.
Then, since $2p> b\ge a$ and since the random variables $Y_n$, $n\ge
1$, are i.i.d., Theorem~1.4.5 of \cite{sam:taqqu} ensures that
$C_1<+\infty$ almost surely.

Let us now focus on $s_{2}$. Observe that
$
\EM{{\llvert  V_n(y) \rrvert }}\le X_n$,
with
\[
X_n=\EM{{\bigl\llvert V_n(x_0) \bigr
\rrvert }}+c_1Y_n
\]
for ${c_1=\sup_{z\in K_{d+1}}\rho(x_0,z)^{\beta} \EM{{\llvert  \log
\EM{{ (1+\rho(x_0,z)^{-1} )}} \rrvert }}^{\eta}}$. Let us
remark that $c_1<+\infty$, by
continuity of $\rho$ on the compact set $\EM{{ \{x_0 \}
}}\times K_{d+1}$.
Moreover, since $V_n(x_0)\in L^{2p}$, $(X_n)_{n\ge1}$ is still a
sequence of i.i.d. variables in $L^{2p}$ and following the same lines
as for equation \eqref{senalpha},
we obtain that, almost surely, for any $x,y\in K_{d+1}$,
\[
s_2(x,y)\le C_2\bigl |\al-\al'\bigr |
\]
with $C_2$ a finite positive random variable.
Let us also note that by equation \eqref{controlrho2},
there exist finite positive constants $c_2$ and $c_3$ such that for any
$x=\EM{{ (\al,u )}}\in K_{d+1}$ and any $ y=\EM{{
(\al',v )}}\in K_{d+1}$,
\[
\EM{{\bigl\llvert \al-\al' \bigr\rrvert }}\le c_2\rho
\EM{{ (x,y )}}^{1/\overline{H}}\le c_3\rho \EM{{ (x,y )}}
\]
since $\overline{H}\le1$.
Hence, since $\be\in(0,1]$,
almost surely, for any $x,y\in K_{d+1}$,
\[
s\EM{{ (x,y )}}\le C \rho\EM{{ (x,y )}}^\beta \log\bigl(1+
\rho(x,y)^{-1}\bigr)^{\max
\EM{{ (\eta,0 )}}}
\]
with $C$ a finite positive random variable. Then assertion 1 follows
from Corollary~\ref{ArgReco}.
\noqed
\end{pf*}

\begin{pf*}{Proof of assertion 2}
Let us choose $p'>0$ such that $1/p'\in(0,1/b-1/\min(2,2p))$. Then,
replacing in the previous lines $s$ by the parameter $r_N$ and Theorem~1.4.5 of \cite{sam:taqqu} by Theorem~2.2 of \cite{COLALE07}, we
obtain: there exists $C$ a finite positive random variable such that
almost surely, for any $x,y\in K_{d+1}$, and for any $N\in\dN$,
\[
r_N\EM{{ (x,y )}}\le C\EM{{ (N+1 )}}^{-1/p'}\rho\EM{{ (x,y
)}}^\beta\log\bigl(1+\rho (x,y)^{-1}\bigr)^{\max\EM{{ (\eta,0 )}}}.
\]
Note also that by Theorem~2.1 in \cite{COLALE07}, almost surely
\[
\sup_{N\in\dN} {N}^{1/p'}\EM{{\bigl\llvert
R_N(x_0) \bigr\rrvert }} <+\infty.
\]
Therefore, by Corollary~\ref{ArgReco2},
almost surely,
\[
\sup_{N\in\dN} N^{1/p'}\sup_{x\in K_{d+1}}\EM{{
\bigl\llvert R_N(x) \bigr\rrvert }}<+\infty,
\]
which concludes the proof.\qed
\noqed
\end{pf*}
\noqed
\end{pf*}

Let us now prove Proposition~\ref{MG}.

\begin{pf*}{Proof of Proposition~\ref{MG}}
Since equation \eqref{controlrho} is fulfilled, there exists $r\in
(0,1)$ such that
$\rho(x,y)\le h_0$ for all $x,y\in K_{d+1}$ with $\EM{{\llVert
x-y\rrVert }}\le r$.
Then, the assumptions done imply that
\[
X_1:= \mathop{\sup_{x,y\in K_{d+1}}}_{0<\EM{{\llVert  x- y\rrVert }}\le
r}
\frac
{\EM{{\llvert  V_1(x)-V_1(y) \rrvert }}}{\rho(x,y)^{\beta}\EM{{\llvert  \log\rho(x,y) \rrvert }}^{\eta
}}\le\sup_{h\in(0,h_0]} \frac{\mathcal{G}(h)}{F(h)
}:=G,
\]
where $F(h):=h^{\beta}\EM{{\llvert  \log h \rrvert }}^{\eta}$. We
assume without loss
of generality that $h_0=2^{-k_0}$ with an integer $k_0\ge1$ is such
that $F$ is increasing on $(0,h_0]$ and equation \eqref{IhTh} holds
for $h\in(0,h_0]$. Then, using the monotonicity of $\mathcal{G}$ and $F$,
\[
G^{2p}\le \sum_{k=k_0}^{+\infty} \sup
_{h\in(2^{-k-1},2^{-k}]} \biggl(\frac
{\mathcal{G}(h)}{F(h)} \biggr)^{2p}\le\max
\bigl(2^\eta,1\bigr)^{2p}\sum_{k={k_0}}^{+\infty}
\biggl(\frac{\mathcal
{G}(2^{-k})}{F(2^{-k})} \biggr)^{2p} .
\]
Therefore,
by equation \eqref{IhTh}
and definition of $F$,
\[
\mathbb{E}\bigl(X_1^{2p}\bigr)\le\mathbb{E}
\bigl(G^{2p}\bigr) \le \max\bigl(2^\eta,1
\bigr)^{2p}\sum_{k={k_0}}^{+\infty} \EM{{
\llvert k\log 2 \rrvert }}^{-1-2p\veps} <+\infty,
\]
which concludes the proof.
\end{pf*}

%s8.3 #&#
\subsection{Proof of Proposition \texorpdfstring{\protect\ref{eq:Lepage}}{4.4}}
\label{eq:Lepage-A}

Let $K_{d+1}=[a,b]\times\prod_{j=1}^d [a_j,b_j]\subset(0,2)\times
\dR^d$.

\begin{pf*}{Proof of assertion 1}
Let us fix an integer $p\ge1$
and consider $x^{(j)}=\EM{{ (\al_j,u^{(j)} )}}\in K_{d+1}$
for each
integer $1\le j\le p$. Then, we set $\vec{x}=\EM{{
(x^{(1)},\ldots, x^{(p)} )}}$. Choosing
$
\mathcal{S}=\EM{{ \{\xi\in\dR^d; m\EM{{ (\xi
)}}>0 \}}}
$
we define $H_{\vec{x}} \dvtx(0,+\infty)\times\mathcal{S}\times\dC
\rightarrow\dC^p$ by
%
%e39 #&#
%
\begin{equation}
\label{Hx} H_{\vec{x}}(r, \xi,g)=\EM{{ \bigl(r^{-1/\al_1}
f_{{\al
_1}}\bigl(u^{(1)},\xi \bigr) m\EM{{ (\xi
)}}^{-1/\al_1} g,\ldots ,r^{-1/\al_p} f_{{\al _p}}
\bigl(u^{(p)},\xi\bigr) m\EM{{ (\xi )}}^{-1/\al_p} g \bigr)}}.
\end{equation}
Let us note that almost surely
\[
\sum_{n=1}^{N} H_{\vec{x}}\EM{{
(T_n,{\xi_n,g_n} )}}=\EM {{
\bigl(S_{m,N}\EM{{ \bigl(x^{(1)} \bigr)}},
\ldots,S_{m}\EM{{ \bigl(x^{(p)} \bigr)}} \bigr)}},
\]
where $S_{m,N}$ is defined by \eqref{SGSN} with $W_n$ given by \eqref
{WnShotnoise}. Then this series
converges almost surely to $\EM{{ (S_{m}\EM{{
(x^{(1)} )}},\ldots,S_{m}\EM{{ (x^{(p)} )}} )}}$.
Since $g_1$ is symmetric,
applying Theorem~2.4 of \cite{Rosinskiu} and using a simple change of
variables ($t=rm(\xi)$) and $\nu\EM{{ (\dR^d\backslash
\mathcal {S} )}}=0$, we obtain that
\[
\forall\la=(\la_1,\ldots,\la_p)\in\dC^p,
\forall z\in\dC ,\qquad \dE\EM{{ \bigl(\mathrm{e}^{\mathrm{i}\Re\EM{{ (\overline{z}\sum
_{j=1}^p \la _jS_m\EM{{ (x^{(j)} )}} )}}} \bigr)}}=\exp
\EM{{ \bigl(I_{\vec{x},\la}\EM{{ (z )}} \bigr)}},
\]
where
\begin{eqnarray*}
&&I_{\vec{x},\la}\EM{{ (z )}}=\int_{(0,+\infty)\times\dR
^d\times\dC} \EM \bigl(
\mathrm{e}^{\mathrm{i}\Re\EM{{ (\overline{z}\langle\la
,J_{\vec{x}}(t,\EM{{ (\xi,g )}})\rangle )}}}
\\
&&\phantom{I_{\vec{x},\la}\EM{{ (z )}}=\int_{(0,+\infty)\times\dR
^d\times\dC} \EM \bigl(}{} - 1-\mathrm{i}\Re\EM{{ \bigl(\overline{z} \bigl\langle
\la,J_{\vec{x}}\bigl(t,\EM {{ (\xi ,g )}}\bigr)\bigr\rangle{
\mathbf{1}}_{\EM{{\llvert  \Re\EM
{{ (\overline{z} \langle\la ,J_{\vec{x}}(t,\EM{{ (\xi
,g )}})\rangle )}} \rrvert }}\le1} \bigr)}} \bigr) \,\mathrm{d}t\nu\EM{{ (\mathrm{d}\xi )}}
\dP_g(\mathrm{d}g)
\end{eqnarray*}
with $\dP_{g}$ the distribution of $g_1$ and
\[
J_{\vec{x}}\bigl(t,\EM{{ (\xi,g )}}\bigr)=\EM{{ \bigl(t^{-1/\al
_1}
f_{{\al _1}}\bigl(u^{(1)},\xi\bigr) g,\ldots,t^{-1/\al_p}
f_{{\al
_p}}\bigl(u^{(p)},\xi\bigr) g \bigr)}}.
\]
Therefore, $I_{\vec{x},\la}$ does not depend on the function $m$, and
then neither does the distribution of the vector $\EM{{ (S_{m}\EM
{{ (x^{(1)} )}},\ldots,S_{m}\EM{{ (x^{(p)}
)}} )}}$.
Since this holds for any $p$ and $\vec{x}$, assertion 1 is established.
\noqed
\end{pf*}

\begin{pf*}{Proof of assertion 2}
Let us now consider the space $B=\cC\EM{{ (K_{d+1},\dC
)}}$ of
complex-valued continuous functions defined on the compact set
$K_{d+1}$. This space is endowed with the topology of the uniform
convergence, so that it is a Banach space.

Let us assume that
$S_{\tilde{m}}$ belongs almost surely to $\mathcal{H}_{\rho}\EM
{{ (K_{d+1}, \beta,\eta )}}\subset B$.
For any $\vec{x}=\EM{{ (x^{(1)},\ldots,x^{(p)} )}}\in
K_{d+1}^p$, in view
of its characteristic function, the vector $\EM  (S_{\tilde
{m}}\EM{{ (x^{(1)} )}},\ldots,\allowbreak  S_{\tilde{m}}\EM{{
(x^{(p)} )}} )$ is infinitely
divisible and its L\'evy measure is given by
\[
%
%\begin{array}{rcl}
F_{\vec{x}}\EM{{ (A )}} = \int_{(0,+\infty)\times\mathcal{S}\times\dC}
{\mathbf {1}}_{A\backslash\EM{{ \{0 \}}}} \EM{{ \bigl({H}_{\vec
{x}}\EM{{ (r,{\xi,g} )}}
\bigr)}} m\EM{{ (\xi )}} \,\mathrm{d}r \nu\EM{{ (\mathrm{d}\xi )}} \dP_g(\mathrm{d}g)
%\end{array}
%
\]
for any Borel set $A\in\cB\EM{{ (\dC^p )}}$. We first
assume that
$
\EM{{ (\al,u )}}\mapsto f_\al(u,\xi)
$
belongs to $B$ for all $\xi\in\dR^d$ so that the function
\begin{eqnarray*}
%
%\begin{array}{rrcl}
H \dvtx (0,+\infty)\times\mathcal{S}\times\dC&
\longrightarrow& B,
\\
(r,{\xi,g}) &\mapsto& \EM{{ \bigl(\EM{{ (\al,u )}}\mapsto r^{-1/\al}
f_\al \EM{{ (u,\xi )}} m\EM {{ (\xi )}}^{-1/\al} \bigr)}}
%\end{array}
%
\end{eqnarray*}
is well-defined.
Since $H_{\vec{x}}$ is defined by \eqref{Hx}, one checks that $\EM
{{ (S_{\tilde{m}}\EM{{ (x )}} )}}_{x\in
K_{d+1}}$ is a $B$-valued infinitely
divisible random variable with L\'evy measure defined by
\[
F\EM{{ (A )}}=\int_{(0,+\infty)\times\mathcal{S}\times
\dC}{\mathbf {1}}_{A\backslash\EM{{ \{0 \}}}} \EM{{
\bigl(H\EM{{ (r,{\xi,g} )}} \bigr)}} m\EM{{ (\xi )}} \,\mathrm{d}r \nu\EM{{ (\mathrm{d}\xi )}}
\dP_g\EM{{ (\mathrm{d}g )}},\qquad A\in\cB(B).
\]
Then, by Theorem~2.4 of \cite{Rosinskiu},
$
\sum_{n=1}^N H\EM{{ (T_n,\EM{{ (\xi_n,g_n )}} )}}
$
converges almost surely in $B$ as $N\to+\infty$. Then, by definition
of $H$, the sequence $\EM{{ (S_{m,N} )}}_{N\in\dN}$
converges in $B$
almost surely. Therefore, its limit $S_{m}$ is almost surely continuous
on $K_{d+1}$.

Let us now consider $\cD\subset K_{d+1}$ a countable dense set in $K_{d+1}$.
Then, since almost surely $S_{\tilde{m}}\in\mathcal{H}_{\rho}\EM
{{ (K_{d+1},\beta,\eta )}}$ and since $S_{m}\stackrel
{\mathrm{fdd}}{=} S_{\tilde{m}}$,
we get that
almost surely
\[
\mathop{\sup_{x,y\in\cD}}_{x\ne y}\frac{\EM{{\llvert
S_{{m}}(x)-S_{{m}}(y) \rrvert }}}{\rho(x,y)^{\beta}\EM{{ [\log
(1+\rho (x,y)^{-1}) ]}}^\eta} <+
\infty.
\]
Then, by continuity of $\rho$, by almost sure continuity of $S_m$ and
by density of $\cD$ on the compact set $K_{d+1}$,
\[
\mathop{\sup_{x,y\in K_{d+1}}}_{x\ne y}\frac{\EM{{\llvert
S_{{m}}(x)-S_{{m}}(y) \rrvert }}}{\rho(x,y)^{\beta}\EM{{ [\log
(1+\rho (x,y)^{-1}) ]}}^\eta} <+
\infty
\]
almost surely, that is, $S_{m}$ belongs almost surely to $\mathcal
{H}_{\rho}\EM{{ (K_{d+1},\beta,\eta )}}$. This
establishes assertion 2
when $\EM{{ (\al,u )}}\mapsto f_{\al}(u,\xi)$ is
continuous for all $\xi
\in\dR^d$.

Assume now that $\EM{{ (\al,u )}}\mapsto f_{\al}(u,\xi)$
is continuous
for $\xi\in\dR^d\backslash\mathcal{N}$ with $\nu\EM{{
(\mathcal {N} )}}=0$ and set
\[
g_{\al}\EM{{ (u,\xi )}}:=f_{\al}(u,\xi){\mathbf{1}}_{\dR
^d\backslash
\mathcal{N}}(
\xi).
\]
Then, almost surely, for all $x=\EM{{ (\al,u )}}\in
(0,2)\times\dR^d$
and all $N\ge1$,
\[
S_{m,N}\EM{{ (x )}}=\sum_{n=1}^N
T_n^{-1/\al} g_{\al}\EM {{ (u,\xi _n
)}}m\EM{{ (\xi_n )}}^{-1/\al} g_n,
\]
and the conclusion follows from the previous lines since $\EM{{
(\al ,u )}}\mapsto g_{\al}\EM{{ (u,\xi )}}$
is continuous on $K_{d+1}$ for all $\xi\in\dR^d$. The proof of
Proposition~\ref{eq:Lepage} is then complete.
\noqed\end{pf*}

%s9 #&#
\section{Applications}
%s9.1 #&#
\subsection{Proof of Proposition \texorpdfstring{\protect\ref{SIOS}}{5.1}}
\label{SIOS-A}

\setcounter{equation}{39}

Let us first note that using Remark~\ref{remComp}, we can and may
assume without loss of generality that $a_1=1$, up to replace $E$ by
$E/a_1$ and $\tau_{{E}}$ by $\tau_{{E/a_1}}^{1/a_1}$.

Let us choose $\zeta>0$ arbitrarily small and consider the Borel
function $\tilde{m}$ defined on $\dR^d$ by
\[
\tilde{m}\EM{{ (\xi )}}= \EM{{\llVert \xi\rrVert }}^{\al_0} {
\mathbf{1}}_{\EM{{\llVert  \xi\rrVert }}\le
A}+ {\tau_{ {{E^t}}}\EM{{ (\xi )}}^{-q(E) } \EM{{
\bigl\llvert \log{\tau_{ {{E^t}}}}(\xi) \bigr\rrvert }}^{-1-\zeta}} {
\mathbf{1}}_{\EM
{{\llVert  \xi\rrVert }}>A}.
\]
Observe that $\tilde{m}$ is positive on $\dR^d\backslash\EM{{
\{0 \}}}$. Then,
\[
0<c=\int_{\dR^d} \tilde{m}\EM{{ (\xi )}} \,\mathrm{d}
\xi=c_{1}+c_{2}
\]
with
\[
c_{1}=\int_{\EM{{\llVert  \xi\rrVert }}\le A}\EM{{\llVert \xi \rrVert
}}^{\al_0} {\,\mathrm{d}\xi} \quad\mbox{and}\quad c_{2}=\int
_{\EM{{\llVert  \xi\rrVert }}> A}{\tau_{ {{E^t}}}\EM {{ (\xi )}}^{-q(E)}
\EM{{\bigl\llvert \log{\tau_{
{{E^t}}}}(\xi) \bigr\rrvert }}^{-1-\zeta}}
{\,\mathrm{d}\xi}.
\]
Let us first observe that $c_{{1}}<\infty$ since $\al_0>0$. To prove
that $c_{2}$ is also a finite constant, we need some tools given in
\cite{thebook,OSSRF}.
As in Chapter~6 of \cite{thebook}, let us consider the norm $\|\cdot\|
_{{E^t}}$ defined by
%
%e40 #&#
%
\begin{equation}
\label{normE} \EM{{\llVert x\rrVert }}_{{E^t}}=\int
_0^1\EM{{\bigl\llVert \theta ^{E^t}x
\bigr\rrVert }}\frac{\mathrm{d}\theta
}{\theta},\qquad\forall x\in{\mathbb R^d}.
\end{equation}
According to the change of variables in polar coordinates (see \cite
{OSSRF}) there exists a finite positive Radon measure $\sigma_{{
E^t}}$ on
$S_{{ E^t}}=\{\xi\in{\mathbb R^d}\dvtx \|\xi\|_{{E^t}}=1\}$ such
that for all
measurable function $\varphi$ non-negative or in $ L^1({\mathbb
R^d},\mathrm{d}\xi)$,
\[
\int_{{\mathbb R^d}}\varphi(\xi)\,\mathrm{d}\xi=\int_0^{+\infty}
\int_{S_{{
E^t}}}\varphi\bigl(r^ {{E^t}}\theta\bigr)
\sigma_{{ E^t}}(\mathrm{d}\theta)r^{q(E)-1}\,\mathrm{d}r.
\]
Applying this change of variables, it follows that $c_{2}<\infty$
since $\zeta>0$.
Hence, $m=\tilde{m}/c$ is well-defined and $\mu\EM{{ (\mathrm{d}\xi
 )}}=m\EM{{ (\xi )}}{\,\mathrm{d}\xi}$ is a probability
measure equivalent to the Lebesgue measure.
Then we may consider $S_m(\alpha_0,u)$ defined by \eqref{SLePage} for
$u\in\dR^d$ so that
$
X_{\al_0}\stackrel{\mathrm{fdd}}{=}d_{\al_0}S_m(\alpha_0,\cdot)
$ with $d_{\al_0}$ given by \eqref{defcste}.

To study the sample path regularity of $S_m(\alpha_0,\cdot)$ on
$K_d=\prod_{j=1}^d[a_j,b_j]$, we apply Proposition~\ref{MG} on
$K_{d+1}=\{\alpha_0\}\times K_d \subset(0,2)\times\dR^d$ for
\[
V_1(\alpha_0,u)=f_{{\al_0}}(u,
\xi_1)m(\xi_1)^{-1/\alpha_0}
\]
with $f_{{\al_0}}$ defined by \eqref{NoyHarmo}.
We recall that here $\xi_1$ is a random vector of ${\mathbb R^d}$ with density
$m$. Therefore let us now check that assumptions of Proposition~\ref
{MG} are fulfilled.

For $h>0$ and $\xi\in{\mathbb R^d}$ we consider
\[
g(h,\xi)=\min\bigl(c_{{ E^t}}\bigl \|h^{E^t}\xi\bigr \|_{{E^t}},1
\bigr)\bigl |\psi_{\al
_0}\EM{{ (\xi )}}\bigr |,
\]
where $c_{{ E^t}}>0$ is chosen such that $|\mathrm{e}^{\mathrm
{i}\langle u,\xi\rangle
}-1|\le c_{{ E^t}}\|\tau_{ {E}}(u)^{E^t}\xi\|_{{E^t}}$.
We consider the quasi-metric defined on $\dR^{d+1}$ by
\[
\rho\bigl((\al,u),\bigl(\al',v\bigr)\bigr)=\bigl |\alpha-
\alpha'\bigr |+\rho_{ {E}}(u,v),\qquad \forall(\alpha,u), \bigl(
\alpha',v\bigr)\in\dR\times{\mathbb R^d},
\]
which clearly satisfies equation \eqref{controlrho}.
By definition of $V_1$, $g$ and $\|\cdot\|_{{E^t}}$, the random field
${\mathcal G}=\EM{{ (g(h,\xi_1) )}}_{h\in[0,+\infty)}$ satisfies
(i) and (ii) of Proposition~\ref{MG}.
It remains to consider assumption (iii).
Let
%
%e41 #&#
%
\begin{equation}
\label{Ih} I(h)=\mathbb{E}\bigl({\mathcal G}(h)^2\bigr)=\int
_{{\mathbb R^d}}g(h,\xi )^2{m}(\xi )^{1-2/\al_0} \,\mathrm{d}\xi.
\end{equation}
Since $\psi_{\al_0}$ satisfies \eqref{CondHS},
\[
I(h)=\int_{\dR^d} m\EM{{ (\xi )}}^{1-2/\al_0} \min\EM{{
\bigl(c_{{E^t}}\EM{{\bigl\llVert h^{E^t}\xi\bigr\rrVert
}}_{{E^t}},1 \bigr)}}^2 \bigl |\psi_{\al_0}(
\xi)\bigr |^2 \,\mathrm{d}\xi\le I_1(h)+I_2(h)
\]
with
\[
I_1(h)= c^{2/\al_0-1}c_{{ E^t}}^2\int
_{\EM{{\llVert  \xi\rrVert
}}\le A} \EM{{\bigl\llVert h^{E^t}\xi\bigr\rrVert
}}_{{E^t}}^2 \EM{{\llVert \xi\rrVert }}^{\al_0-2} \bigl |
\psi_{\al_0}(\xi )\bigr |^2 {\,\mathrm{d}\xi},
\]
where $A$ is given by the condition \eqref{CondHS}, and
\[
I_2(h)=c^{2/\al_0-1}c_\psi\int_{\dR^d}
\min\EM{{ \bigl(c_{{ E^t}}\EM{{\bigl\llVert h^{E^t}\xi\bigr
\rrVert }}_{{E^t}},2 \bigr)}}^2 \tau_{ {{E^t}}}\EM{{ (\xi
)}}^{-q(E)-2\beta} \EM{{\bigl\llvert \log\tau_{ {{E^t}}}(\xi) \bigr\rrvert
}}^{(1+\zeta)(2/\al
_0-1)}\,\mathrm{d}\xi.
\]
From Lemma~3.2 of \cite{BL09} there exists {a finite constant} $C_1>0$
such that for all $h\in(0,\mathrm{e}^{-1}]$
\[
I_1(h)\le C_1 h^{2a_1}\bigl |\log(h)\bigr |^{2(d-1)}.
\]
Moreover, using again the change of variables in polar coordinates,
there exists {a finite constant} $C_2>0$ such that for all $h\in(0,\mathrm{e}^{-1}]$,
\[
I_2(h)\le C_2 h^{2\beta}\bigl |\log(h)\bigr |^{(1+\zeta)(2/\al_0-1)}.
\]
Since $\beta<a_1$, one find {a finite constant} $C_3>0$ such that
%
%e42 #&#
%
\begin{equation}
\label{controlIh} I(h)\le C_3h^{2\beta}\bigl |\log(h)\bigr |^{2(1+\zeta)(1/\al_0-1/2)}.
\end{equation}
Hence, assumption (iii) of Proposition~\ref{MG} is also
fulfilled and applying this proposition, it follows that \eqref{UBV}
is satisfied with $\beta$ and $\eta=1/\alpha_0+\varepsilon$, for
all $\varepsilon>0$.
Then, by Theorem~\ref{CrMom},
almost surely $S_m\in{\mathcal H}_{\rho}(K_{d+1},\beta,1/\al
_0+1/2+\varepsilon)$. By definition of $\rho$ and $K_{d+1}$, this
means that
\[
S_{m}\EM{{ (\al_0,\cdot )}}\in{\mathcal
H}_{\rho_{{E}}}(K_{d},\beta ,1/\al_0+1/2+\varepsilon).
\]
In particular, $S_{m}\EM{{ (\al_0,\cdot )}}$ is
continuous on $K_d$.
Then, since
$d_{\al_0} S_m(\al_0,\cdot)\stackrel{\mathrm{fdd}}{=}X_{\al_0}$, $X_{\al
_0}$ is stochastically continuous and almost surely
\[
C:= \sup_{u,v\in\mathcal{D}, u\neq v}\frac{|X_{\al
_0}(u)-X_{\al_0}(v)|}{\tau_{ {E}}(u-v)^\beta\EM{{ [\log
(1+\tau_{ {E}}(u-v)^{-1}) ]}}^{1/\al_0+1/2+\varepsilon
}}<+\infty,
\]
where $\mathcal{D}\subset K_d$ is a countable dense in $K_d=\prod_{j=1}^d[a_j,b_j]$. So let us write $\Omega^*$ this event and let us
define a modification of $X_{\alpha_0}$ on $K_d$.

First, if $\omega\notin\Omega^*$, we set $X^*_{\al_0}\EM{{
(u )}}\EM{{ (\omega )}}=0$ for all $u\in K_d$. Let
us now fix $\omega\in\Omega^*$.
Then, we set
\[
X^*_{\al_0}\EM{{ (u )}}\EM{{ (\omega )}}=X_{\al_0}\EM{{ (u )}}
\EM{{ (\omega )}},\qquad \forall u\in\mathcal{D}.
\]
Let us now consider $u\in K_d$. Then, there exists
$u^{(n)}\in\cD$ such that
$\lim_{n\to+\infty}u^{(n)}=u$. It follows that,
\begin{eqnarray*}
&&\EM{{\bigl\llvert X^*_{\al_0}\EM{{ \bigl(u^{(n)} \bigr)}}(
\omega )-X^*_{\al_0}\EM{{ \bigl(u^{(m)} \bigr)}}(\omega) \bigr
\rrvert }}
\\
&&\quad\le C(\omega) \tau_{ {E}}\bigl(u^{(n)}-u^{(m)}
\bigr)^\beta \EM{{ \bigl[\log\bigl(1+\tau_{ {E}}
\bigl(u^{(n)}-u^{(m)}\bigr)^{-1}\bigr)
\bigr]}}^{1/\al
_0+1/2+\varepsilon},
\end{eqnarray*}
so that $\EM{{ (X^*_{\al_0}\EM{{ (u^{(n)} )}}\EM
{{ (\omega )}} )}}_n$ is a Cauchy
sequence and hence converges. We set
\[
X_{\al_0}^*\EM{{ (u )}}\EM{{ (\omega )}}=\lim_{n\to+\infty}X_{\al
_0}^*
\EM{{ \bigl(u^{(n)} \bigr)}}\EM{{ (\omega )}}.
\]
Remark that this limit does not depend on the choice of $\EM{{
(u^{(n)} )}}_n$ and that
$X_{\al_0}^*\EM{{ (\cdot )}}(\omega)$ is then well-defined on $K_d$.
Observe also that, by stochastic continuity of $X_{\al_0}$, $X_{\al
_0}^*$ is a modification of $X_{\al_0}$.
Moreover, by continuity of $\tau_{ {E}}$,
\[
C(\omega)= \sup_{u,v\in K_d, u\neq v}\frac{|X_{\al
_0}^*(u)(\omega)-X_{\al_0}^*(v)(\omega)|}{\tau_{ {E}}(u-v)^\beta
\EM{{ [\log(1+\tau_{ {E}}(u-v)^{-1}) ]}}^{1/\al
_0+1/2+\varepsilon
}}<+\infty
\]
for all $\omega\in\Omega$
and $X_{\al_0}^*$ is continuous on $K_d$.
This concludes the proof.%\qed
%$\square$

%s9.2 #&#
\subsection{Multistable random fields}
\label{SeMStable-A}

This section is devoted to the proofs of the results stated in
Section~\ref{SeMStable}. Let us first establish Proposition~\ref{ModCMS}.

\begin{pf*}{Proof of Proposition~\ref{ModCMS}}
Since $\tilde{\rho
}$ satisfies equation \eqref{controlrho}, so does $\rho$. Then,
assumptions of Theorem~\ref{CrMom} are fulfilled, which implies that
$\EM{{ (S_{m,N} )}}_{N\in\dN}$ converges uniformly to
$S_m$ on
$K_{d+1}=[a,b]\times K_d$. Therefore, $\EM{{ (\tilde
{S}_{m,N} )}}_{N\in
\dN}$ converges uniformly to $\tilde{S}_m$ on $K_d$ since $ \tilde
{S}_{m,N}(u)=S_{m,N}\EM{{ (\al(u),u )}}$ and $ \tilde
{S}_m(u)=S_m\EM{{ (\al(u),u )}}$ and $\al$ is continuous.

Moreover, by Theorem~\ref{CrMom}
there exists a finite positive random variable $C$ such that for any
$u,v\in K_d$,
\[
\EM{{\bigl\llvert \tilde{S}_m\EM{{ (u )}}-\tilde{S}_m
\EM {{ (v )}} \bigr\rrvert }}\le C \rho\EM{{ \bigl(x(u),x(v) \bigr)}}^\beta
\EM{{ \bigl[\log\EM{{ \bigl(1+\rho\EM{{ \bigl(x(u),x(v) \bigr)}}^{-1}
\bigr)}} \bigr]}}^{\max(\eta,0)+1/2},
\]
where $x(w)=\EM{{ (\al(w),w )}}$.
Moreover, by definition of $\rho$ and since $\al\in\mathcal
{H}_{\tilde{\rho}}\EM{{ (K_d,1,0 )}}$, there exists a
finite positive
constant $c_1$ such that
\[
\forall u,v\in K_d,\qquad\rho\EM{{ \bigl(x(u),x(v) \bigr)}}\le
c_1\tilde{\rho}\EM{{ (u,v )}}.
\]
Let us now recall that since $\tilde{\rho}$ is continuous on the
compact set $K_d\times K_d$,
$
M=\sup_{u,v\in K_d} \tilde{\rho}\EM (u,\allowbreak  v )<+\infty$.
Then, up to change $C$, for all $u,v\in K_d$,
\[
\EM{{\bigl\llvert \tilde{S}_m\EM{{ (u )}}-\tilde{S}_m
\EM {{ (v )}} \bigr\rrvert }}\le C \tilde{\rho }\EM{{ (u,v )}}^\beta
\EM{{ \bigl[\log\EM{{ \bigl(1+\tilde{\rho}\EM{{ (u,v )}}^{-1} \bigr)}}
\bigr]}}^{\max\EM{{ (\eta,0 )}}+1/2}
\]
since $h\mapsto h^{\beta} \log(1+h^{-1})^{\max\EM{{ (\eta
,0 )}}+1/2}$
is increasing around $0$ and bounded on $[0,M]$. Assertion~1 is then
proved. Moreover, assertion 2 is a direct consequence of assertion 2
of Theorem~\ref{CrMom}. The proof is then complete.
\end{pf*}

Let us conclude this paper by the proof of Corollary~\ref{MSV}.

\begin{pf*}{Proof of Corollary~\ref{MSV}}
Let $K_d=\prod_{j=1}^d
[a_j,b_j]\subset\dR^d$ and $u_0\in K_d$. Let us set
\[
a=\min_{{K_{ d}}}\al,\qquad b=\max_{{K_{ d}}} \al
\quad\mbox {and}\quad K_{d+1}=[a,b]\times K_d \subset(0,2)
\times\dR^d.
\]
Let us first note that
Assumption~\ref{HLep} is fulfilled
with $K_1=[a,b]$ and then $\tilde{S}_m$ is well-defined. Let us now
consider $\rho_{ {E}}$ and $\tau_{ {E}}$ as defined in Example~\ref
{tauE}. Then we set
\[
\tilde{m}\EM{{ (\xi )}}=\frac{c_{\zeta}}{\tau_{
{{E^t}}}\EM{{ (\xi  )}}^{{q(E)}} \EM{{\llvert  \log\tau
_{ {{E^t}}}\EM{{ (\xi )}} \rrvert }}^{1+\zeta}},
\]
with $\zeta>0$ a parameter chosen arbitrarily small. Therefore, let us consider
\[
\tilde{V}_n\EM{{ (\al,u )}}= f_{\al}\EM{{ (u,\tilde{
\xi}_n )}} \tilde {m}\EM{{ (\tilde{\xi}_n
)}}^{-1/\al},
\]
where $(\tilde{\xi_n})_{n\ge1}$ is a sequence of i.i.d. random
variables with common distribution $\tilde{\mu}\EM{{ (\mathrm{d}\xi
 )}}=\tilde
{m}(\xi)\,\mathrm{d}\xi$. The sequence $(\tilde{\xi_n})_{n\ge1}$ is assumed
to be independent from $\EM{{ (T_n,g_n )}}_{n\ge1}$. Then,
Assumption~\ref
{HSN} is fulfilled. Moreover,
\[
\EM{{\bigl\llvert \tilde{V}_n\EM{{ (\al,u )}}-
\tilde{V}_n\EM {{ \bigl(\al',v \bigr)}} \bigr\rrvert
}}\le\EM{{\bigl\llvert \tilde{V}_n\EM {{ (\al,u )}}-
\tilde{V}_n\EM{{ (\al,v )}} \bigr\rrvert }}+\EM{{\bigl\llvert \tilde
{V}_n\EM{{ (\al,v )}}-\tilde{V}_n\EM{{ \bigl(
\al',v \bigr)}} \bigr\rrvert }}.
\]
Let us set
\[
C_1= \mathop{\sup_{u,v\in K_d}}_{0<\EM{{\llVert  u- v\rrVert }}\le
r}\sup
_{\al\in[a,b]} \frac{\EM{{\llvert  \tilde{V}_1\EM{{ (\al,u )}}-\tilde
{V}_1\EM{{ (\al ,v )}} \rrvert }}}{{\rho_{{E}}}\EM
{{ (u,v )}}\EM{{\llvert  \log{\rho_{E}\EM{{
(u,v )}}} \rrvert }}^{{\eta}}}
\]
and
\[
C_2= \mathop{\sup_{\al,\al'\in[a,b]} }_{\al\ne\al'}\sup
_{u\in K_d}\frac{\EM{{\llvert  \tilde{V}_1\EM
{{ (\al,u )}}-\tilde {V}_1\EM{{ (\al',u )}}
\rrvert }}}{\EM{{\llvert  \al-\al' \rrvert }} },
\]
where $r>0$ and the choice of $\eta\in\dR$ is given below.
Then, for any $x=\EM{{ (\al,u )}}\in K_{d+1}$ and any
$y=\EM{{ (\al ',v )}}\in K_{d+1}$ such that $\EM{{\llVert  x-y\rrVert }}\le r$,
\begin{eqnarray*}
%
%\begin{array}{rcl}
\EM{{\bigl\llvert \tilde{V}_1\EM{{ (x )}}-
\tilde{V}_1\EM {{ (y )}} \bigr\rrvert }}&\le& \EM{{
(C_1+C_2 )}}\EM {{ \bigl({\rho}_{{E}}\EM{{
(u,v )}} \EM{{\bigl\llvert \log {\rho}_{{E}}\EM{{ (u,v )}} \bigr\rrvert
}}^{\eta}+ \EM {{\bigl\llvert \al-\al' \bigr\rrvert }}
\bigr)}}
\\
&\le& c_1\EM{{ (C_1+C_2 )}} \rho\EM{{ (x,y
)}}\EM{{\bigl\llvert \log{\rho}\EM{{ (x,y )}} \bigr\rrvert }}^{\eta},
%\end{array}
%
\end{eqnarray*}
where $c_1\in(0,+\infty)$ is a finite constant and $\rho\EM{{
(x,y )}}={\rho}_{{E}}\EM{{ (u,v )}} +\EM{{\llvert
\al-\al' \rrvert }}$. Then, to apply
assertion 1 of Proposition~\ref{ModCMS} with $\tilde{\rho}={\rho
}_{{E}}$ and $\beta=1$, it suffices to establish that $C_1,C_2$ and
$ \tilde{V}_1(\al(u_0),u_0)\in L^{2}$ (since $b< 2$).

Let us first deal with $\tilde{V}_1\EM{{ (\al\EM{{
(u_0 )}},u_0 )}}$.
Using polar coordinates associated with $E^t$ (see \cite{thebook}),
\[
\dE\EM{{ \bigl(\EM{{\bigl\llvert \tilde{V}_1\EM{{ \bigl(\al\EM{{
(u_0 )}},u_0 \bigr)}} \bigr\rrvert }}^2
\bigr)}}\le c_2 \int_{0}^{+\infty} \min
\EM{{ \bigl(\EM{{\bigl\llVert t^{E^t}\bigr\rrVert }},1
\bigr)}}^2t^{-3}\EM{{\llvert \log t \rrvert
}}^{2(1+\zeta)/\al
\EM{{ (u_0 )}}-1} \,\mathrm{d}t
\]
with $c_2\in(0,+\infty)$.
Hence,
Lemma~2.1 of \cite{OSSRF} proves that $V_1\EM{{ (\al\EM{{
(u_0 )}},u_0 )}}\in
L^2$ for any choice of~$\zeta$.

Let us now consider the random variable $C_1$. By homogeneity and
continuity of $\psi$, there exists a finite positive constant $c_3$
such that for any $u,v\in K_d$,
\[
\sup_{\al\in[a,b]} \EM{{\bigl\llvert \tilde{V}_1\EM{{ (
\al ,u )}}-\tilde{V}_1\EM{{ (\al,v )}} \bigr\rrvert }}\le
c_3 \EM{{\bigl\llvert \mathrm{e}^{\mathrm{i}\langle u-v,\tilde{\xi}_1\rangle }-1 \bigr\rrvert }}
Z_1
\]
with
\[
Z_1=\tau_{{ {E^t}}}\EM{{ (\tilde{\xi}_1
)}}^{-1} \max \EM{{ \bigl(\EM{{\bigl\llvert \log \tau_{{ {E^t}}}
\EM{{ (\tilde {\xi}_1 )}} \bigr\rrvert }}^{(1+\zeta)/a}, \EM{{\bigl
\llvert \log \tau _{{ {E^t}}}\EM{{ (\tilde{\xi}_1 )}} \bigr
\rrvert }}^{(1+\zeta
)/b} \bigr)}}.
\]
Combining the proofs of Propositions~\ref{MG} and~\ref{SIOS}, we
obtain that for any $\veps>0$, choosing $r$ small enough,
\[
\dE\EM{{ \biggl(\EM{{ \biggl[\mathop{\sup_{u,v\in K_d}}_{0<\EM
{{\llVert  u- v\rrVert }}\le r}
\frac{\EM{{\llvert  {\mathrm
{e}^{\mathrm{i}\langle u-v,\tilde{\xi}_1\rangle}-1} \rrvert }} Z_1}{{\rho_{{E}}}\EM{{ (u,v )}}\EM{{\llvert  \log{\rho_{E}\EM
{{ (u,v )}}} \rrvert }}^{1/a+\veps} } \biggr]}}^2 \biggr)}}<+\infty.
\]
This implies that for any $\veps>0$, $C_1\in L^2$ for $\eta=1/a+\veps
$ and $\zeta$ well-chosen.

Let us now study $C_2$. Since $K_d$ is a compact set, using polar
coordinates and the Mean Value theorem, we have
\[
\sup_{v\in K_d}\EM{{\bigl\llvert \tilde{V}_1\EM{{ (
\al,v )}}-\tilde{V}_1\EM{{ \bigl(\al ',v \bigr)}} \bigr
\rrvert }}\le c_4 \EM {{\bigl\llvert \al-\al' \bigr
\rrvert }} Z_2
\]
with $Z_2=\min\EM{{ (\EM{{\llVert  \tau_{{E^t}}\EM{{
(\xi_n )}}^{E^t}\rrVert }},1 )}} Z_1
\EM{{\llvert  \log\tau_{{E^t}}\EM{{ (\xi_n )}}+c_5
\rrvert }}$ and $c_4$ and $c_5$ two
finite positive constants. Using polar coordinates, one checks that
$Z_2\in L^2$, which implies that $C_2\in L^2$.

Therefore, for any $\veps>0$, assumptions of assertion 1 of
Proposition~\ref{ModCMS} are fulfilled for a well-chosen $\zeta$.
This implies that almost surely, for any $\veps>0$, $\tilde{S}_m\in
\mathcal{H}_{\rho_{E}}\EM{{ (K_d,1,1/a+1/\veps )}}$
with $a=\min_{K_d} \al$.
Hence, for any $\veps>0$, $\tilde{S}_m\in\mathcal{H}_{\rho_{E},B(u_0,r)}\EM{{ (u_0,1,1/\al(u_0)+1/2+1/\veps )}}$ for
$r$ small enough.
This concludes the proof.
\end{pf*}
\end{appendix}

%\section*{Acknowledgements}

%\begin{supplement}%[id=suppA]
%\slink[doi]{10.3150/00-BEJXXXXSUPP} %[doi,text={...}] - jei reikia
%suskaldyti doi
%\sfilename{BEJ000\_supp.pdf}
%\end{supplement}

%\bibitem[\protect\citeauthoryear{}{()}]{r1}
%\bibitem{r1}

% imsref loaded by audrone.aklyte, 2014-06-10 13:31:30
% imsref loaded by audrone.aklyte, 2014-06-10 13:33:37
%

\printhistory
\end{document}